\documentclass[12pt]{article}
\usepackage{amssymb}
\usepackage{amsmath}
\usepackage{amsthm}
\usepackage{enumitem}
\usepackage{graphicx}
\begin{document}

\newtheorem{theorem}{Theorem}
\newtheorem{conjecture}[theorem]{Conjecture}
\newtheorem{lemma}[theorem]{Lemma}

\title{Dynamics of Simple Balancing Models with
State Dependent Switching Control.}
\author{
D.J.W.~Simpson, R.~Kuske and Y.-X.~Li\\
Department of Mathematics\\
University of British Columbia\\
Vancouver, BC, V6T1Z2\\
Canada}
\maketitle

\begin{abstract}
Time-delayed control in a balancing problem may be a nonsmooth 
function for a variety of reasons.
In this paper we study a simple model of the control of an
inverted pendulum by either a connected movable cart or an applied torque
for which the control is turned off when the pendulum
is located within certain regions of phase space.
Without applying a small angle approximation for
deviations about the vertical position, we see structurally
stable periodic orbits which may be attracting or repelling.
Due to the nonsmooth nature of the control,
these periodic orbits are born in various discontinuity-induced bifurcations.
Also we show that a coincidence of switching events
can produce complicated periodic and aperiodic solutions.
\end{abstract}

%=====================================================================
\section{Introduction}
\label{sec:INTRO}

The subject of balance has received considerable recent attention.
Local measurements of muscles \cite{LoLa02b,CaMo05,LoMa05b}
and experiments that highlight the influence of vision
\cite{DaSt93,MaPo03}
have led to an improved understanding 
of the key physiological elements in human balancing tasks.
From a theoretical perspective,
progress has been made in analyzing systems with time-delay \cite{Er09}.
Time-delay is used to model the reaction time of the controlling mechanism
and is a near ubiquitous element of mathematical models of balancing tasks.
A current challenge is to incorporate new experimental
observations into mathematical models and interpret the results.
Time-delayed balance control is also a fundamental problem in
robotics \cite{BaMi10,Br08,BuLe05},
however the control strategies used in engineering
are typically distinct from those identified in physiology \cite{MiCa09}.

Time-delayed PD control
(for which the applied forces are determined from measurements of position (P)
and its derivative (D))
has been used in models of the control of a stiff beam or rod that is
attached at its base to a controllable cart \cite{LaCa05,SiKr04}
and stick balancing \cite{StKo00}.
In both cases for small time-delay it is possible to choose
control parameters so that the model
is successfully directed to the vertical position.
With an increase in the delay time, for fixed control parameters,
the vertical position becomes unstable which may result in
stable oscillations about the vertical position.
Beyond a critical value of the delay, the control is incapable
of stabilizing the system at the vertical position.
A similar response has been found with other smooth control laws \cite{SiKr05}.
Some authors have used more complex models to incorporate
additional physical features such as friction \cite{CaCr08}.
In general, time-delayed control models are inherently difficult
to analyze because they are infinite-dimensional.
However, dynamics near local bifurcations of delay differential equations
are well described by low-dimensional systems of ordinary differential equations
derived via a centre manifold analysis \cite{Ca09,SiKr04b}.

Recently several researchers have proposed nonsmooth balancing models.
Experimental observations of human balance in quiet standing \cite{MiTo09}
and studies of human balancing tasks \cite{MiOh09,CaMi02}
suggest a switching control due to intermittent muscle movements \cite{LoGa06}.
Specifically, experiments suggest a state-dependent control
such as the ``drift and act'' method of \cite{MiOh09}
which simply turns the control off when the system nears equilibrium.
Muscle control may be active or passive,
the latter refers to muscles which intrinsically resist motion
away from equilibrium. %\cite{LoLa02b,CaMo05}.
The presence of control is potentially detrimental to obtaining equilibrium
if the system is naturally approaching equilibrium;
in \cite{AsTa09} the authors consider a switching control
that lessens this effect.
In \cite{StIn06}, the authors present a control that
acts only after waiting for a time longer than the delay in order
to gain sufficient information from the system to be able
to perform more effective control.
Hysteretic control laws, which are common in temperature control,
have been considered in balancing models \cite{Si06,WaCh08}.
In mechanical systems small gaps between gears create
backlash which another source of nonsmoothness
and in the absence of friction eliminates the possibility of perfect stable equilibrium
at the vertical position \cite{KoSt00}.
From an engineering or robotics viewpoint,
a switching control may require less cost and be easier to implement
mechanically or may be necessary due to a nonzero sampling time.
An additional benefit is that nonsmooth control laws are often able to
stabilize the system for arbitrarily large delay in simple mathematical models.
For an introduction to switching in control systems, see for instance \cite{Li03}.

Mathematical models that incorporate both time-delay and switching conditions
are usually particularly difficult to analyze, yet there are common mechanisms
that induce a transition from simple to complex dynamics
and are characteristic of such systems.
For instance a periodic orbit may develop a tangency with a switching condition,
or undergo switching at times that differ by exactly the delay time of the system
\cite{Si06,CoDi07,SiKo10b}.
Both scenarios correspond to a codimension-one bifurcation of the periodic orbit.

The purpose of this paper is to investigate the effect of
the application of control with two different
switching rules on simple balancing models.
We study nonlinear equations of motion for a pendulum
combined with time-delayed PD control by a force applied
either by a movable cart or as a torque.
We consider two different ON/OFF switching rules for the application of the control.
For the majority of this paper we analyze
a switching rule that involves both position and velocity,
a later section of the paper is devoted to a switching rule based on position only.
The goal is to reveal and understand novel dynamics
resulting from the combined effects of time-delay,
switching and nonlinearity which are generic and therefore
expected to be prevalent in a wide range of balancing systems.

The remainder of this paper is organized as follows.
The equations of motion investigated are stated in \S\ref{sub:EQNS}.
In \S\ref{sub:SWITCH} we introduce two switching rules
that divide phase space into various ON and OFF regions.
We then focus on the first switching rule.
Section \ref{sub:NO} summarizes dynamics when the delay time is zero.
In this case orbits may become ``stuck'' to a switching manifold.
In the presence of small time-delay this {\em sliding motion}
becomes rapid switching motion, which we refer to as {\em zigzag} motion
about the switching manifold, \S\ref{sub:BASIC}.
This motion corresponds to a jittery motion
of the pendulum restricted to one side of the vertical position
maintained by an intermittent application of the control.
The time-delay may alternatively induce {\em spiral} motion
corresponding to oscillations about the vertical position,
but for short time-delay zigzag dynamics dominate \S\ref{sub:ZIGZAG}.
The bifurcation structure is detailed in \S\ref{sub:BIFSET}.
In particular we prove that both stable and unstable zigzag periodic orbits
may be born in discontinuity-induced bifurcations.
These bifurcations are described formally through asymptotic expansions
in \S\ref{sub:EXPAN}.
Homoclinic zigzag orbits are the subject of \S\ref{sub:HOMOC}.

In \S\ref{sec:LARGE} we consider longer values of time-delay
for which numerics reveal complex dynamics that are
not explained by the small delay asymptotic expansions.
We describe a novel bursting-like attractor
that exhibits different behaviours on two distinct time-scales, \S\ref{sub:BURST}.
In \S\ref{sub:FOUR} we map out the stability of the pendulum
at the vertical position in the plane of the control parameters by linearizing the
equations of motion.

The switching rule based on position only is studied in \S\ref{sec:OO}.
This rule corresponds to the controller neglecting
control when near the vertical position.
We again compute asymptotic expansions and prove the existence of
stable periodic orbits.
Finally a summarizing discussion is presented in \S\ref{sec:CONC}.

%=====================================================================
\section{Simple Balancing Models}
\label{sec:MODEL}

Here we detail the mathematical models studied in this paper.
To give our results wide applicability we have chosen
to study dimensionless, inverted pendulum-based equations of motion
that have been used as models in both human balancing tasks
and mechanical systems.
For simplicity we consider all motions to be restricted to a plane
and ignore both friction and noise.
Section \ref{sub:EQNS} introduces the equations of motion
and PD control,
\S\ref{sub:SWITCH} details two ON/OFF control mechanisms
and describes basic dynamics.
Dynamics in the absence of delay are described in \S\ref{sub:NO}.

%---------------------------------------------------------------------
\subsection{Equations of Motion}
\label{sub:EQNS}

The task of vertically balancing a long stick
has been modelled as an inverted pendulum
with an applied torque
\cite{MiCa09,StIn06,StKo00}.
Upon an appropriate time scaling and other reductions
the equation of motion may be written simply as
\begin{equation}
\ddot{\theta} - \sin(\theta) = F \;,
\label{eq:e1}
\end{equation}
where $\theta$ is a dimensionless quantity representing
the angular displacement of the stick from vertical
and $F$ denotes pivot control due to, say,
the finger or hand of the human performing the balancing actions.

Equation (\ref{eq:e1}) also provides a simple
model of human postural sway, where $F$ represents
ankle torque \cite{AsTa09}.
For postural sway it is important to treat
non-control components of $F$.
In particular, ankle torque has an intrinsic passive stiffness
which provides some stability but is regarded as inadequate
to maintain quiet standing \cite{LoLa02b,CaMo05}.
The ankle joint also provides damping.
If human postural sway is overdamped it may be
suitable to model the motion with a first order differential
equation
\cite{MiOh09,EuMi96}.
Analyses of (\ref{eq:e1}) provide a basis for more complex motions
such as balancing with hip movements and bipedal walking
in humans and robots.

Planar dynamics of a vertical rod controlled
by a moving cart have been modelled by
\begin{equation}
\left( 1 - \frac{3m}{4} \cos^2(\theta) \right) \ddot{\theta} +
\frac{3m}{8} \dot{\theta}^2 \sin(2\theta) - \sin(\theta) = F \cos(\theta) \;,
\label{eq:cart}
\end{equation}
\cite{LaCa05,SiKr04},
where $m = \frac{M_{\rm pendulum}}{M_{\rm pendulum} + M_{\rm cart}}$
denotes the fraction of the mass of the system that belongs to the pendulum
and friction is ignored.
If the cart is much more massive than the pendulum,
i.e.~$M_{\rm cart} \gg M_{\rm pendulum}$, then $m=0$ is a useful approximation
and the equation of motion simplifies to
\begin{equation}
\ddot{\theta} - \sin(\theta) = F \cos(\theta) \;.
\label{eq:e2}
\end{equation}
Equation (\ref{eq:cart}) exhibits dynamics similar to
(\ref{eq:e2}) for small values of $m$.
%The opposite limiting case, $m=1$, provides less simplification.
A cart model with friction is studied in \cite{CaCr08}.

In this paper we study both (\ref{eq:e1}) and (\ref{eq:e2}).
For convenience we let $\phi = \dot{\theta}$ and
write
\begin{equation}
\begin{split}
\dot{\theta} &= \phi \;, \\
\dot{\phi} &= \sin(\theta) + F G(\theta) \;,
\end{split}
\label{eq:e}
\end{equation}
where $G(\theta) = 1$ or $G(\theta) = \cos(\theta)$.
For physical scenarios %such as human postural sway during quiet standing
that involve only small changes in angular displacement,
it is suitable to linearize equation of motions in $\theta$.
We do not use a linear (small angle) approximation in $\theta$ in order to
investigate invariant solutions created in bifurcations at $\theta = 0$.

The control force, $F$, is a time-delayed function of $\theta$ and $\phi$.
In the context of human balancing tasks time-delay models neural transmission time.
One of the simplest control laws is PD control \cite{Du05,BuLe05}:
\begin{equation}
F_{\rm ON} = a \theta(t-\tau) + b \phi(t-\tau) \;,
\label{eq:PD}
\end{equation}
where $a$ and $b$ are scalar control parameters and $\tau \ge 0$ is
the delay time.
For inverted pendulum-type problems with PD control that is continuously applied,
the vertical position is typically stable for control parameters
that lie in a roughly D-shaped region in the $(a,b)$-plane \cite{SiKr04,StKo00}.
Pitchfork bifurcations and Hopf bifurcations form the boundary of this region.
Outside the region there may exist stable oscillations about the vertical position,
complex dynamics or even failure for the system to attain a
physically meaningful bounded solution.
The area of the D-shaped region reduces as the value of $\tau$ is increased.
In \cite{SiKr04b},
the authors study (\ref{eq:cart}) with (\ref{eq:PD})
and find that the D-shaped region shrinks to a point at some critical
value of $\tau$.

%---------------------------------------------------------------------
\subsection{ON/OFF control}
\label{sub:SWITCH}

By ON/OFF control we mean simply that
at some times the control is implemented,
whereas at other times the control is absent, i.e.~$F=0$.
In this paper we study two different 
ON/OFF controls based on position in the $(\theta,\phi)$-plane.
In contrast the ON/OFF mechanism in the so-called ``act-and-wait''
strategy \cite{StIn06,InKo10,StIn07} is based on the time elapsed.
When control is applied we use the PD control (\ref{eq:PD}).
We refer to (\ref{eq:e}) with $F = 0$ as the OFF system 
and (\ref{eq:e}) with (\ref{eq:PD}) as the ON system.

\begin{figure}[b!]
\begin{center}
\includegraphics[width=7.2cm,height=6cm]{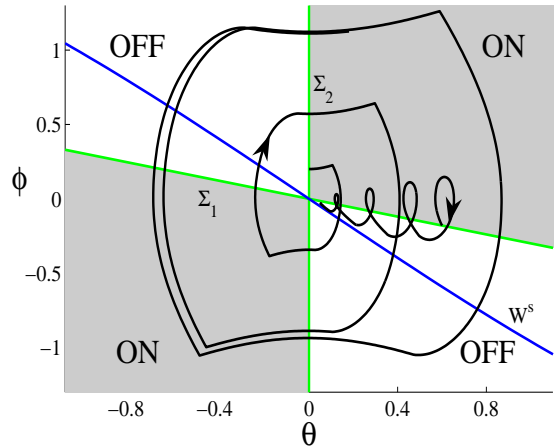}
\caption{
The $(\theta,\phi)$-plane for (\ref{eq:e})-(\ref{eq:PD}) with (\ref{eq:c1}).
Two trajectories were numerically computed when
$s = -0.3$, $\tau = 0.5$, $(a,b) = (1.5,4)$ and $G(\theta) = \cos(\theta)$.
One trajectory zigzags about $\Sigma_1$ and tends to the origin,
the other trajectory spirals out from the origin
limiting upon a stable periodic orbit.
$W^s$ is the stable manifold of the origin for the OFF system.
\label{fig:onOffc1}
}
\end{center}
\end{figure}

The two control laws we analyze are 
\begin{eqnarray}
F &=& \left\{ \begin{array}{lc}
F_{\rm ON} \;, & \theta(t-\tau) \big( \phi(t-\tau) - s \theta(t-\tau) \big) > 0 \\
0 \;, & {\rm otherwise}
\end{array} \right. \;, \label{eq:c1} \\
F &=& \left\{ \begin{array}{lc}
F_{\rm ON} \;, & |\theta(t-\tau)| > \sigma \\
0 \;, & {\rm otherwise}
\end{array} \right. \;, \label{eq:c2}
\end{eqnarray}
where $s \le 0$ and $\sigma > 0$.
Equations (\ref{eq:c1}) and (\ref{eq:c2})
are respectively taken from \cite{AsTa09} and \cite{MiOh09,KoGl10},
which focus on human postural sway.
The control law (\ref{eq:c1}) defines two switching manifolds
\begin{equation}
\begin{split}
\Sigma_1 &= \{ (\theta, s \theta) ~|~ \theta \in \mathbb{R} \} \;, \\
\Sigma_2 &= \{ (0,\phi) ~|~ \phi \in \mathbb{R} \} \;,
\end{split}
\label{eq:swManC1}
\end{equation}
that divide the $(\theta,\phi)$-plane into four regions.
We refer to these regions as ON and OFF regions,
as indicated in Fig.~\ref{fig:onOffc1}.
Similarly (\ref{eq:c2}) defines two switching manifolds
\begin{equation}
\begin{split}
\Sigma_3 &= \{ (\sigma, \phi) ~|~ \phi \in \mathbb{R} \} \;, \\
\Sigma_4 &= \{ (-\sigma, \phi) ~|~ \phi \in \mathbb{R} \} \;,
\end{split}
\label{eq:swManC2}
\end{equation}
that divide the $(\theta,\phi)$-plane into an OFF region
that is a vertical strip centred about $\theta = 0$,
and two ON regions, Fig.~\ref{fig:onOffc2}.
The system (\ref{eq:e})-(\ref{eq:PD})
retains the usual symmetry by $(\theta,\phi) \mapsto -(\theta,\phi)$
with either switching condition.

%%%%%%%%%%%%%%%%%%%%%%%%%%%%%%%%%%%%%%%%%%%%%%%%%%%%%%%%%%%%%%%%%%%%%%
%%%%%%%%%%%%%%%%%%%%%%%%%%%%%%%%%%%%%%%%%%%%%%%%%%%%%%%%%%%%%%%%%%%%%%

The system (\ref{eq:e})-(\ref{eq:PD}) with either (\ref{eq:c1}) or (\ref{eq:c2})
is a piecewise-smooth, discontinuous, delay differential equation system.
Due to the presence of time-delay the phase space of the system
is infinite-dimensional \cite{Er09,HaLu93,DiVa95},
nevertheless it is convenient to picture dynamics in the $(\theta,\phi)$-plane.
It is difficult to consider all possible initial conditions
(which are curves, $(\theta(t),\phi(t))$ for $t \in [-\tau,0]$).
However, when $F=0$, (\ref{eq:e})
reduces to a two-dimensional ODE system (the OFF system).
Consequently, if the point $(\theta(0),\phi(0))$ lies in an OFF region
and the trajectory $(\theta(t),\phi(t))$ has been in this OFF region
for a time equal to at least $\tau$,
then $F=0$ for all $t \in [0,\tau]$
and so the trajectory at any positive time is independent of
its location at any time prior to $t=0$.
Hence, with this assumption, the initial condition
may be thought of as merely the location of the trajectory at $t=0$.
Since it is straight-forward to understand the dynamics of the OFF system,
we consider as initial conditions
only points on switching manifolds at which the vector field
of the OFF system points into the neighbouring ON region.
Admittedly this restriction omits some behaviour of the system,
but we believe it captures all the important and physically meaningful dynamics.

%---------------------------------------------------------------------
\subsection{Dynamics in the Absence of Delay}
\label{sub:NO}

Here we describe the behaviour of (\ref{eq:e})-(\ref{eq:c1}) when $\tau = 0$.
We begin with brief analyses of the individual OFF and ON systems when $\tau = 0$,
then describe the effects of switching rule (\ref{eq:c1}).

%~~~~~~~~~~~~~~~~~~~~~~~~~~~~~~~~~~~~~~~~~~~~~~~~~~~~~~~~~~~~~~~~~~~~~
\subsubsection*{The OFF system}

The OFF system, given by (\ref{eq:e}) with $F=0$,
represents a classical inverted pendulum and is independent of the time-delay.
The system is Hamiltonian and the function
\begin{equation}
H(\theta,\phi) = \frac{1}{2} \phi^2 + \cos(\theta) \;,
\label{eq:hamiltonian}
\end{equation}
is a suitable Hamiltonian function for this system
(i.e.~$H(\theta(t),\phi(t))$ is constant for any solution).
In \S\ref{sub:EXPAN} we will use (\ref{eq:hamiltonian})
to measure the variation of trajectories in the presence of 
small delay over the course of individual zigzag oscillations.
The equilibria of the OFF system are $(n\pi,0)$, for $n \in \mathbb{Z}$.
If $n$ is even, the equilibrium is a saddle;
if $n$ is odd, the equilibrium is a centre.

%~~~~~~~~~~~~~~~~~~~~~~~~~~~~~~~~~~~~~~~~~~~~~~~~~~~~~~~~~~~~~~~~~~~~~
\subsubsection*{The ON system}

In the absence of delay, the ON system for (\ref{eq:e})-(\ref{eq:PD}) is
\begin{equation}
\begin{split}
\dot{\theta} &= \phi \;, \\
\dot{\phi} &= \sin(\theta) - (a \theta + b \phi) G(\theta) \;,
\end{split}
\label{eq:onNo}
\end{equation}
which corresponds to instantaneous PD control.
For both choices of $G$, the origin is an equilibrium of (12) and 
it can be characterized  with a standard stability calculation.
For $a<1$ the origin is a saddle,
otherwise it is a node or a focus.
The node or focus is stable [unstable] for $b>0$ [$b<0$].

For $G(\theta) = \cos(\theta)$, (\ref{eq:onNo}) has infinitely many equilibria,
and for $a>1$ the equilibrium $(\pm \theta_{\rm cos}^*,0)$ is a saddle, where
$\theta_{\rm cos}^*$ denotes the smallest
positive $\theta$-value for the equilibria.
When $b > 0$, the value $a=1$ corresponds to
a subcritical pitchfork bifurcation at the origin.

For $G(\theta) = 1$, the origin is the only equilibrium of (\ref{eq:onNo}) 
for $a \ge 1$. 
For $0 \le a<1$ the equilibrium $(\pm \theta_1^*,0)$
is stable when $b > 0$, where $\theta_1^*$
denotes the smallest positive $\theta$-value for the equilibria.
Consequently, when $b > 0$, the value $a = 1$ corresponds to
a supercritical pitchfork bifurcation at the origin.
For both choices of $G$ the pitchfork bifurcation
forms the left boundary of the $D$-shaped stability region, \S\ref{sub:EQNS}.

%~~~~~~~~~~~~~~~~~~~~~~~~~~~~~~~~~~~~~~~~~~~~~~~~~~~~~~~~~~~~~~~~~~~~~
\subsubsection*{The full system}

Here we analyze (\ref{eq:e})-(\ref{eq:c1}) when $\tau = 0$, which may be written as
\begin{equation}
\left[ \begin{array}{c}
\dot{\theta} \\ \dot{\phi} 
\end{array} \right] = \left\{ \begin{array}{lc}
\left[ \begin{array}{c}
\phi \\ \sin(\theta) - (a\theta+b\phi) G(\theta)
\end{array} \right] \;, & \theta (\phi - s\theta) > 0 \\
\left[ \begin{array}{c}
\phi \\ \sin(\theta)
\end{array} \right] \;, & {\rm otherwise}
\end{array} \right. \;.
\label{eq:noC1}
\end{equation}
This system is a Filippov system \cite{Fi64,Fi88}
because it is discontinuous on the switching manifolds,
$\Sigma_1$ and $\Sigma_2$ (\ref{eq:swManC1}).
Dynamical behaviour that lies entirely within either an ON or OFF region
is determined purely by either the ON or OFF system.
Thus here we focus on trajectories that impact a switching manifold.
For impacts on $\Sigma_2$,
since $\dot{\theta} = \phi$ for both ON and OFF systems,
trajectories simply arrive at $\Sigma_2$ from the neighbouring OFF region
and immediately enter the adjacent ON region.

For impacts on $\Sigma_1$, we identify regions
where the trajectory remains on $\Sigma_1$ for some time. 
A section of $\Sigma_1$ along which the
ON vector field points into the OFF region and the
OFF vector field points into the ON region 
is known as an {\em attracting sliding region} \cite{Fi88,DiBu08}.
A trajectory that arrives at an attracting
sliding region becomes stuck on the switching manifold and {\em slides}.
The manner by which sliding dynamics evolve
is usually defined by Filippov's method, \cite{Fi64,Fi88,DiBu08,LeVa00},
which we explain below.

First let us locate sliding regions on $\Sigma_1$.
At an end of a sliding region the vector field of either the OFF system
or the ON system is tangent to $\Sigma_1$.
Such a point is referred to as a {\em grazing point}
and a trajectory that exhibits this tangency is a {\em grazing trajectory}.
The slope of the vector field of the OFF system on $\Sigma_1$, 
$\frac{\dot{\phi}}{\dot{\theta}} = \frac{\sin(\theta)}{s\theta}$,
 is tangent to $\Sigma_1$ when
\begin{equation}
\mathcal{F}(\theta) \equiv \frac{\sin(\theta)}{\theta} - s^2 = 0 \;.
\label{eq:thOFFgrazC1}
\end{equation}
Similarly, the slope of the vector field of the ON system on $\Sigma_1$ is
tangent to $\Sigma_1$ when
\begin{equation}
\mathcal{G}(\theta) \equiv \frac{\sin(\theta)}{\theta} - s^2 - (a+bs) G(\theta) = 0 \;.
\label{eq:thONgrazC1}
\end{equation}
Roots of $\mathcal{F}$ and $\mathcal{G}$
correspond to possible grazing points and boundaries of sliding regions.
For any $s \in (-1,0]$, $\mathcal{F}$
has a unique root, $\theta_{\rm graz}^{\rm OFF} \in (0,\pi]$, and
the OFF vector field points into the ON region on
$\Sigma_1$ for $0 < \theta < \theta_{\rm graz}^{\rm OFF}$.
However this point does not influence physically meaningful dynamics
unless the value of $s$ is close to $-1$ because
whenever $s > -\sqrt{\frac{2}{\pi}} \approx -0.7979$,
$\theta_{\rm graz}^{\rm OFF} > \frac{\pi}{2}$.

To identify roots of $\mathcal{G}$, (\ref{eq:thONgrazC1}),
we consider the two cases of $G(\theta)$ separately and
omit some messy but elementary calculations.
For $G(\theta) = 1$, 
there is a unique root for $\frac{2}{\pi} < a + bs + s^2 < 1$ given by
$\theta_{\rm graz}^{\rm ON} \in (0,\frac{\pi}{2})$.
It follows that the subset of $\Sigma_1$ for which
$\theta_{\rm graz}^{\rm ON} < \theta < \frac{\pi}{2}$
is an attracting sliding region.
%as shown in Fig.~\ref{fig:noAllC1}-B.
For $G(\theta) = \cos(\theta)$, when $a > 1-bs-s^2$
sliding occurs on $\Sigma_1$ for $0 < \theta < \theta_{\rm graz}^{\rm ON}$
where $\theta_{\rm graz}^{\rm ON}$ is the smallest positive root of $\mathcal{G}$.

Dynamics on an attracting sliding region are governed by
the unique convex combination of 
the ON and OFF vector fields that is tangent to the region at each point.
We write
\begin{equation}
\left[ \begin{array}{c} \dot{\theta} \\ \dot{\phi} \end{array} \right]_{\rm slide} =
(1-q) \left[ \begin{array}{c} \dot{\theta} \\ \dot{\phi} \end{array} \right]_{\rm OFF} +
q \left[ \begin{array}{c} \dot{\theta} \\ \dot{\phi} \end{array} \right]_{\rm ON} \;,
\label{eq:FilippovMethod}
\end{equation}
where $\left[ \begin{array}{c} \dot{\theta} \\ \dot{\phi} \end{array} \right]_{\rm OFF}$
and $\left[ \begin{array}{c} \dot{\theta} \\ \dot{\phi} \end{array} \right]_{\rm ON}$
refer to (\ref{eq:e}) with $F=0$ and (\ref{eq:onNo}) respectively,
and $q$ is a $\theta$-dependent scalar quantity determined by the requirement:
${\dot{\phi}_{\rm slide}}=s{\dot{\theta}_{\rm slide}}$.
Upon substituting $\phi = s \theta$ into (\ref{eq:FilippovMethod}),
${\dot{\phi}_{\rm slide}}=s{\dot{\theta}_{\rm slide}}$ yields
$q = \frac{\sin(\theta) - s^2 \theta}{(a+bs)\theta\cos(\theta)}$
which leads to the explicit solution
\begin{equation}
\left[ \begin{array}{c} \theta_{\rm slide}(t) \\ \phi_{\rm slide}(t) \end{array} \right] =
\theta_0 {\rm e}^{st}
\left[ \begin{array}{c} 1 \\ s \end{array} \right] \;.
\nonumber
\end{equation}
Therefore when $s < 0$ attracting sliding trajectories approach the origin;
when $s=0$, attracting sliding regions are intervals of equilibria.

%%%%%%%%%%%%%%%%%%%%%%%%%%%%%%%%%%%%%%%%%%%%%%%%%%%%%%%%%%%%%%%%%%%%%%
\begin{figure}[b!]
\begin{center}
\setlength{\unitlength}{1cm}
\begin{picture}(10,12.9)
\put(0,6.7){\includegraphics[width=10cm,height=6.2cm]{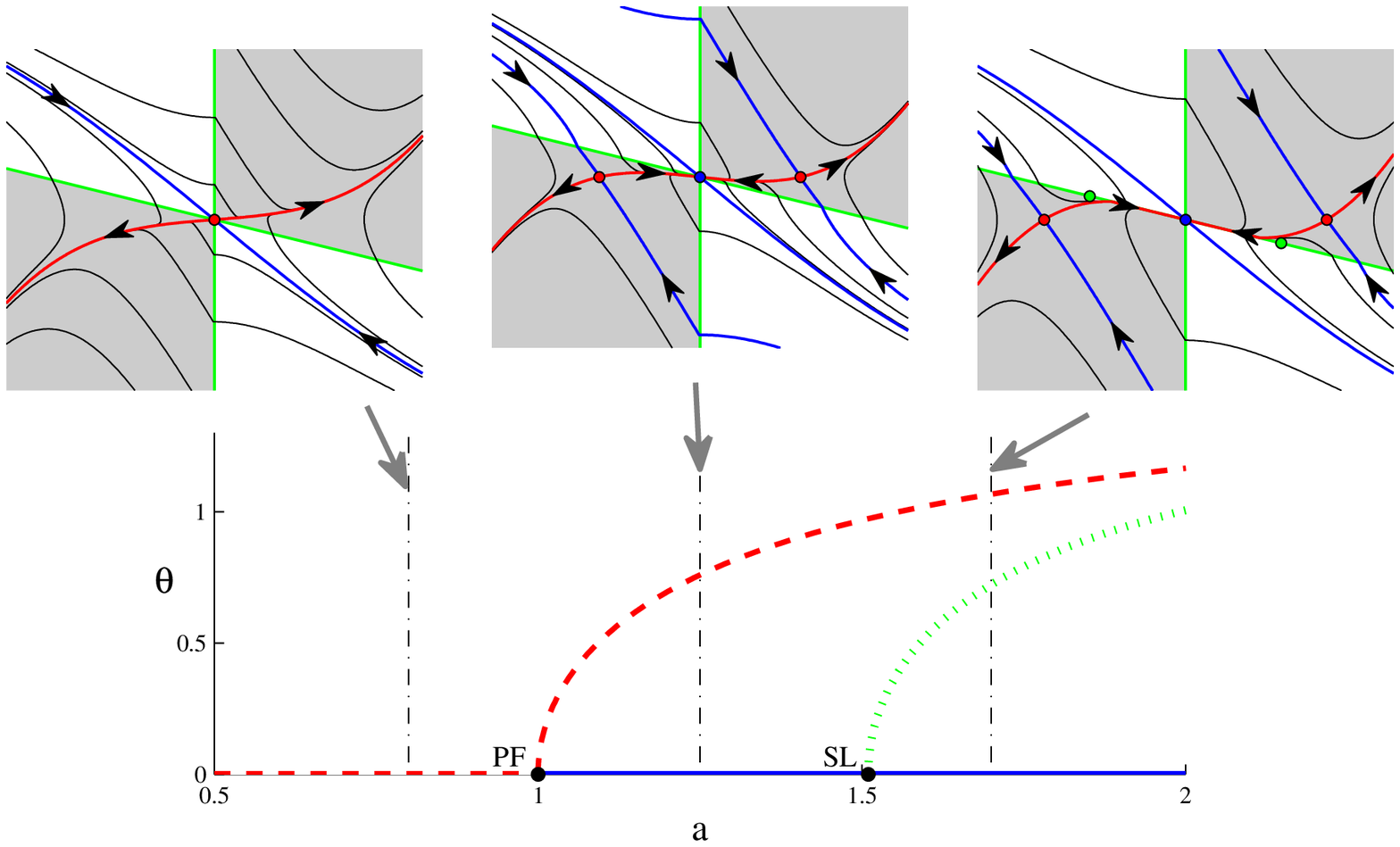}}
\put(0,0){\includegraphics[width=10cm,height=6.2cm]{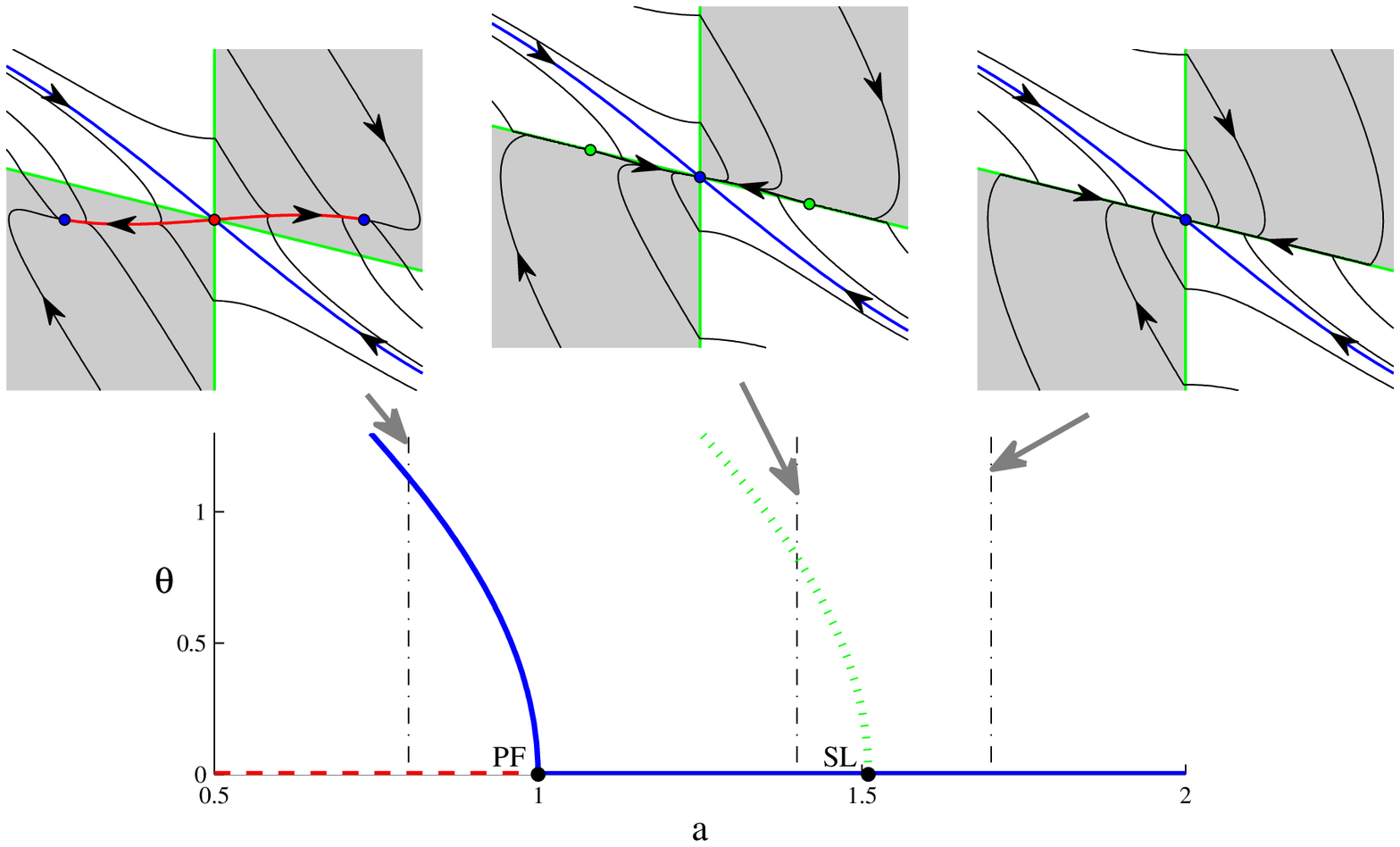}}
\put(.7,12.9){\large \sf \bfseries A}
\put(.7,6.2){\large \sf \bfseries B}
\end{picture}
\caption{
Bifurcation diagrams of the system (\ref{eq:e})-(\ref{eq:c1})
in the absence of delay, (\ref{eq:noC1}), with $s=-0.3$ and $b=2$
for $G(\theta) = \cos(\theta)$ in panel A and $G(\theta) = 1$ in panel B.
Solid and dashed curves denote stable and unstable equilibria, respectively.
The dotted curves indicate $\theta_{\rm graz}^{\rm ON}$,
which is a root of $\mathcal{G}(\theta)$ (\ref{eq:thONgrazC1})
and corresponds to the onset of sliding.
PF - pitchfork-like bifurcation; SL - grazing point intersects origin.
Included are representative phase portraits
with dash-dot lines indicating the corresponding values of $a$.
The system exhibits dynamics similar to that shown here
for a significant range of $s$ and $b$ values;
some bounds on these values are given in the text.
\label{fig:noAllC1}
}
\end{center}
\end{figure}
%%%%%%%%%%%%%%%%%%%%%%%%%%%%%%%%%%%%%%%%%%%%%%%%%%%%%%%%%%%%%%%%%%%%%%

Fig.~\ref{fig:noAllC1}-A illustrates typical
dynamics of (\ref{eq:noC1})
when $G(\theta) = \cos(\theta)$.
There are two bifurcations.
At $a=1$, two saddle equilibria are created that exist to the right of this line.
This bifurcation is a pitchfork bifurcation of the ON system but
not a {\em bona fide} pitchfork bifurcation of the full system
because the origin is a non-differentiable point.
Throughout the paper we refer to $a=1$ as a
{\em pitchfork-like bifurcation}.
At $a=1-bs-s^2$ (labeled SL in the figure),
a sliding region is created on $\Sigma_1$
that grows in size with increasing $a$.
Panel B summarizes dynamics when $G(\theta) = 1$.

%=====================================================================
\section{Dynamics with Delay}
\label{sec:DELAY}

%---------------------------------------------------------------------
\subsection{Zigzag and Spiral Dynamics}
\label{sub:BASIC}

As discussed in \S\ref{sub:SWITCH},
we consider the forward evolution of a point 
on $\Sigma_1$ or $\Sigma_2$ in order to identify the basic behavior
of (\ref{eq:e})-(\ref{eq:c1}). 
Typically the corresponding trajectory immediately enters
the neighbouring ON region, and by symmetry 
we may assume that this is the ON region with $\theta > 0$.
Notice $\dot{\theta} > 0$ whenever $\phi > 0$,
thus the trajectory cannot exit the ON region through $\Sigma_2$.

One possibility is that the trajectory enters the ON region and
remains in the ON region for all time.
In this case a physically meaningful stable solution is not attained.
Roughly speaking this occurs for small values of the control parameters.
Alternatively the trajectory intersects $\Sigma_1$
for the first time at some $t_1 > 0$.
In this case generically the trajectory then resides
in the OFF region, with $\theta > 0$, for some nonzero time.
Assuming the trajectory intersects either $\Sigma_1$ or $\Sigma_2$
at a later time, let $t_2$ denote the
earliest such intersection time.

If $t_2 \ge t_1 + \tau$, 
that is, if the trajectory has been in the OFF region for
a time greater than or equal to the delay time, $\tau$,
then the fate of the trajectory at times later than $t_2$
is independent of its location at any time prior to $t_2$.
Therefore in this case the location of the trajectory at $t = t_2$
may be thought of as a new initial point, \S\ref{sub:SWITCH}.
If instead $t_2 < t_1 + \tau$,
the behaviour of the trajectory is dependent on the control
for $t \in (t_2,t_1+\tau)$.
Trajectories with this property require more effort 
to analyze beyond $t = t_1 + \tau$ and may be particularly complicated but occur
commonly only when $\tau$ is relatively large, see \S\ref{sub:BURST}.

Here we consider the former case, $t_2 \ge t_1 + \tau$, 
and classify two basic types of delay-induced dynamics,
as noted qualitatively in \cite{AsTa09}.
At the time $t = t_1 + \tau$ the applied control is switched off.
It is instructive to consider the location of the trajectory at this time,
i.e.~the point $(\theta(t_1+\tau),\phi(t_1+\tau))$,
in relation to the stable manifold of the origin for the OFF system,
$W^s$, shown in Fig.~\ref{fig:onOffc1}.
If the point $(\theta(t_1+\tau),\phi(t_1+\tau))$ lies above $W^s$
then the point $(\theta(t_2),\phi(t_2))$ lies on $\Sigma_1$;
if $(\theta(t_1+\tau),\phi(t_1+\tau))$ lies below $W^s$
then $(\theta(t_2),\phi(t_2))$ lies on $\Sigma_2$.
(In the special case that the point $(\theta(t_1+\tau),\phi(t_1+\tau))$ lies on $W^s$
($H(\theta(t_1+\tau),\phi(t_1+\tau)) = 1$),
the trajectory coincides with $W^s$ after $t_1$ and never exits the OFF region.)
If $(\theta(0),\phi(0))$ and $(\theta(t_2),\phi(t_2))$ both lie on $\Sigma_1$
we refer to the part of the trajectory between these two points
as a {\em zigzag oscillation}.
Similarly if $(\theta(0),\phi(0))$ and $(\theta(t_2),\phi(t_2))$ both lie on $\Sigma_2$
we refer to the same part of the trajectory 
as half a {\em spiral oscillation}.
Zigzag oscillations often come in succession, as do spiral oscillations.
We have not observed sustained switching between zigzag and spiral motion.
Zigzag motion is typical for small values of $\tau$
and spiral motion is typical for large values of $\tau$,
however, for a carefully tuned combination of the control parameters
zigzag and spiral motion may coexist, as in Fig.~\ref{fig:onOffc1},

%---------------------------------------------------------------------
\subsection{Domination of zigzag trajectories for small delay}
\label{sub:ZIGZAG}

\begin{figure}[h!]
\begin{center}
\setlength{\unitlength}{1cm}
\begin{picture}(8,4.3)
\put(0,0){\includegraphics[width=8cm,height=4.3cm]{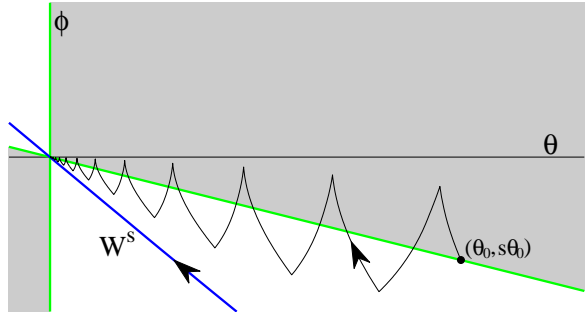}}
\end{picture}
\caption{
The forward orbit of a point on $\Sigma_1$
for (\ref{eq:e})-(\ref{eq:c1}) with $s = -0.3$, $\tau = 0.25$
and $(a,b) = (2.5,2)$.
\label{fig:alwaysZigzag}
}
\end{center}
\end{figure}

Let us consider the forward orbit of a point on $\Sigma_1$, $(\theta_0,s\theta_0)$,
as shown in Fig.~\ref{fig:alwaysZigzag}, for the system (\ref{eq:e})-(\ref{eq:c1}).
With a small enough time-delay ($\tau < -s$ is sufficient),
the orbit, as governed by the OFF system,
does not cross the $\theta$-axis before the control is applied.
If the applied control is sufficiently large so that $\dot{\phi}(t) < 0$ in (\ref{eq:e}),
the orbit abruptly changes heading and soon reintersects $\Sigma_1$.
After a time $\tau$ beyond this reintersection,
the control is switched off.
As long as the applied control is not so strong that
the orbit has overshot $W^s$,
the orbit then continues back to $\Sigma_1$ and the process repeats.
For all time $\theta(t)$ is strictly decreasing because $\phi(t)$
is always negative.
In other words the orbit zigzags into the origin.

Numerical investigations suggest that as long as $|s|$ is not too large,
say $-0.4 < s < 0$, and $\tau < -s$,
then there exists a range of parameters for which
the forward orbit of any point $(\theta_0,s\theta_0)$,
with $|\theta_0| < \theta_{\rm cos}^*$ for $G(\theta) = \cos(\theta)$
and $|\theta_0| < \frac{\pi}{2}$ for $G(\theta) = 1$,
zigzags into the origin as in Fig.~\ref{fig:alwaysZigzag}.
With stronger control (specifically larger values of $a$)
orbits may cross $W^s$ producing spiral dynamics.
With larger delay the orbits enter the first quadrant of the $(\theta,\phi)$-plane.
In this case orbits may still zigzag, but not necessarily approach
the origin.
The next section investigates these dynamics.

%%%%%%%%%%%%%%%%%%%%%%%%%%%%%%%%%%%%%%%%%%%%%%%%%%%%%%%%%%%%%
%%%%%%%%%%%%%%%%%%%%%%%%%%%%%%%%%%%%%%%%%%%%%%%%%%%%%%%%%%%%%

%---------------------------------------------------------------------
\subsection{Bifurcation sets}
\label{sub:BIFSET}

In the previous section we argued that if $\tau < -s$,
then for appropriately chosen $a$ and $b$ orbits simply zigzag into the origin.
For the remainder of this section we consider small values of $s$
so that the condition $\tau < -s$ may not be satisfied.

\begin{figure}[b!]
\begin{center}
\setlength{\unitlength}{1cm}
\begin{picture}(12,12)
\put(0,0){\includegraphics[width=12cm,height=12cm]{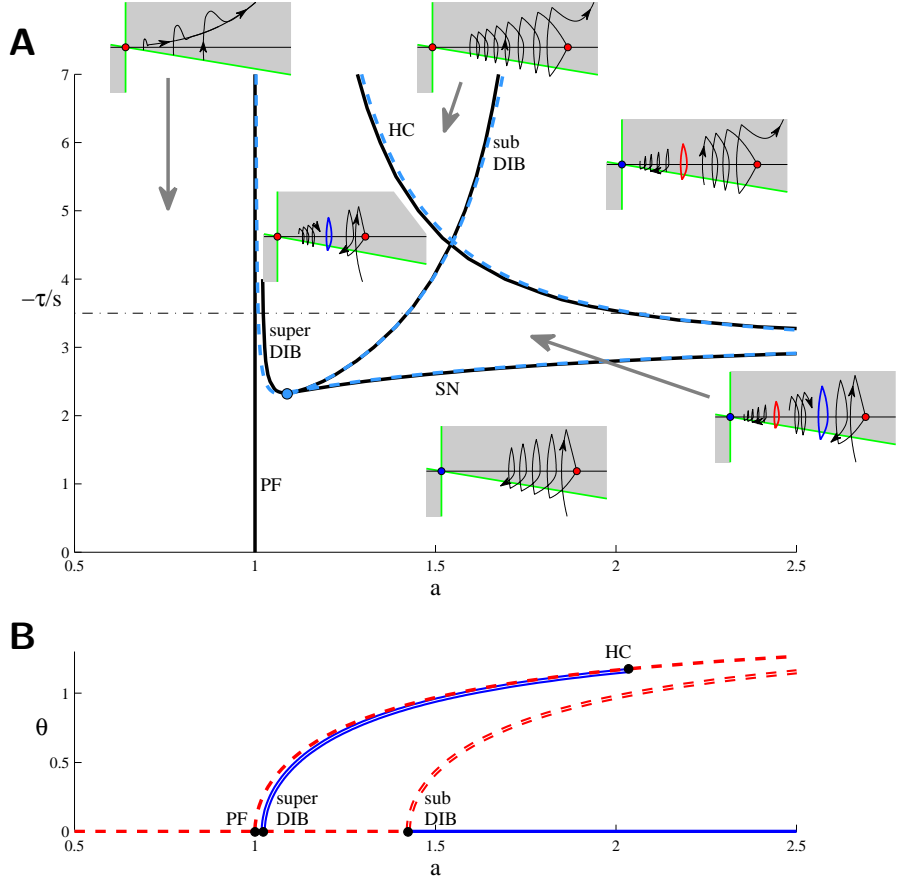}}
\put(.2,11.3){\large \sf \bfseries A}
\put(.2,3.4){\large \sf \bfseries B}
\end{picture}
\caption{
Panel A is a bifurcation set of (\ref{eq:e})-(\ref{eq:c1})
when $s=-0.01$, $b=2$ and $G(\theta) = \cos(\theta)$.
The solid curves are the result of numerical computations.
The dashed curves correspond to equations derived in \S\ref{sub:EXPAN}
and \S\ref{sub:HOMOC}.
PF - pitchfork-like bifurcation;
SN - saddle-node bifurcation of periodic orbits;
HC - homoclinic bifurcation;
DIB - discontinuity-induced bifurcation
({\em super} and {\em sub} are abbreviations for supercritical
and subcritical, respectively).
Included are six representative sketches of trajectories
in the $(\theta,\phi)$-plane.
The sketches are exaggerated for clarity -
in reality $\theta(t)$ changes by an extremely small amount
over each zigzag.
Panel B is a bifurcation diagram corresponding to the
horizontal dash-dot line in panel A (i.e.~for $\tau = 0.035$).
The solid [dashed] curves denote stable [unstable] equilibria.
The double curves correspond to the maximum $\theta$-values
of periodic orbits with line style indicating stability in the same fashion.
\label{fig:smallAllE2C1}
}
\end{center}
\end{figure}

Fig.~\ref{fig:smallAllE2C1}-A is a numerically computed bifurcation set of
(\ref{eq:e})-(\ref{eq:c1}) when $G(\theta) = \cos(\theta)$.
For this figure we have fixed the value of $s$ at $-0.01$;
by setting the vertical axis to $\frac{\tau}{-s}$
the picture is roughly unchanged for different small $s < 0$.
We fixed $b=2$ for Fig.~\ref{fig:smallAllE2C1} and
numerically have observed that the bifurcation structure
is qualitatively the same for different values of $b>0$.
The figure shows dynamics only for $\theta > 0$;
identical dynamics occurs for $\theta < 0$ since the system is symmetric.

%%%%%%%%%%%%%%%%%%%%%%%%%%%%%%%%%%%%%%%%%%%%%%%%%%%%%%%%%%%%%
%%%%%%%%%%%%%%%%%%%%%%%%%%%%%%%%%%%%%%%%%%%%%%%%%%%%%%%%%%%%%

For sufficiently small $\tau$, and $a>1$,
Fig.~\ref{fig:smallAllE2C1}-A reflects the
results of \S\ref{sub:ZIGZAG}: orbits zigzag about $\Sigma_1$
and approach the origin.
An increase in $\tau$ leads to the creation of
a pair of zigzag periodic orbits in a
classical saddle-node bifurcation.
One periodic orbit is stable, the other is unstable.
If $a \gtrsim 1.54$,
a further increase in $\tau$ destroys the stable zigzag periodic orbit
via a homoclinic connection to the saddle equilibrium, $(\theta_{\rm cos}^*,0)$.
Otherwise an increase in $\tau$ destroys the unstable zigzag periodic orbit
by a collision of the periodic orbit with the origin.
This is a discontinuity-induced bifurcation
(labelled sub DIB in Fig.~\ref{fig:smallAllE2C1}-A)
that in some ways resembles a subcritical Hopf bifurcation.
Indeed the amplitude of the periodic orbit grows
at a rate proportional to the square root of parameter change
(shown by Fig.~\ref{fig:smallAllE2C1}-B)
but it does not correspond to the occurrence of
purely imaginary stability multipliers,
nor does the periodic orbit encircle the equilibrium.
Furthermore, the system has identical dynamics for $\theta < 0$,
so a second unstable zigzag periodic orbit is created simultaneously
and exists for $\theta < 0$.

The curve of discontinuity-induced bifurcations
and the curve of saddle-node bifurcations are tangent to one another
at their point of intersection.
The criticality of the discontinuity-induced bifurcation changes here.
That is, along the discontinuity-induced bifurcation curve for $a \lesssim 1.09$,
a stable zigzag periodic orbit existing for $\theta > 0$ is created.
Beyond $-\frac{\tau}{s} = 4$ we were unable to numerically continue the
bifurcation curve when $s = -0.01$.
Certainly as $\tau$ increases the discontinuity-induced bifurcation
approaches the pitchfork-like bifurcation ($a=1$).
A more detailed investigation of this area of parameter space is
presented in \S\ref{sub:BURST} for $s = -0.1$ and suggests that
complex dynamics may occur here.

Numerical simulations indicate that for different values of $b>0$,
the system exhibits a basic bifurcation structure identical
to that shown in Fig.~\ref{fig:smallAllE2C1}.
As $b \to 0^+$ the intersection of the HB and SN curves
appears to approach $(a,-\frac{\tau}{s}) = (1,2)$
and the HC curve seems to limit on the SN curve and $a=1$.
For the values of $b$ and $s$ corresponding to
Fig.~\ref{fig:smallAllE2C1}, when $\tau = 0$,
sliding occurs for $a > 1-bs-s^2 = 1.0199$
which is so close to $a=1$ that we have chosen not to
indicate it in the figure.

%%%%%%%%%%%%%%%%%%%%%%%%%%%%%%%%%%%%%%%%%%%%%%%%%%%%%%%%%%%%%
\begin{figure}[b!]
\begin{center}
\setlength{\unitlength}{1cm}
\begin{picture}(12,12)
\put(0,0){\includegraphics[width=12cm,height=12cm]{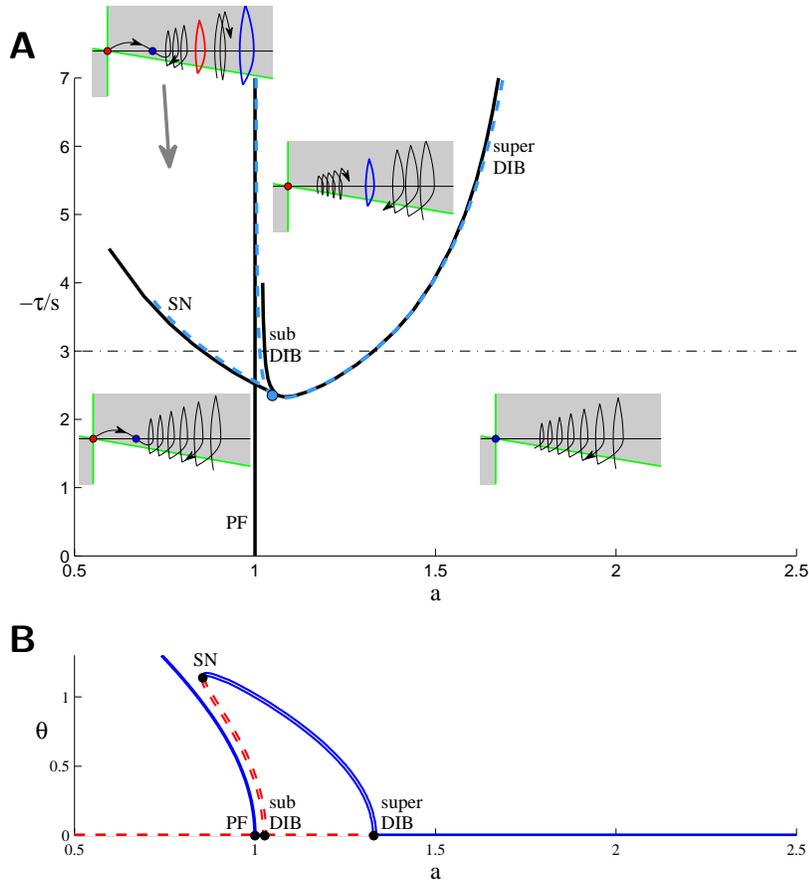}}
\put(.2,11.3){\large \sf \bfseries A}
\put(.2,3.4){\large \sf \bfseries B}
\end{picture}
\caption{
A bifurcation set and bifurcation diagram of (\ref{eq:e})-(\ref{eq:c1})
when $s=-0.01$, $b=2$ and $G(\theta) = 1$.
The meanings of the abbreviations and the significance of the line styles are the
same as in Fig.~\ref{fig:smallAllE2C1}.
\label{fig:smallAllE1C1}
}
\end{center}
\end{figure}
%%%%%%%%%%%%%%%%%%%%%%%%%%%%%%%%%%%%%%%%%%%%%%%%%%%%%%%%%%%%%

A bifurcation set for the same parameter values as Fig.~\ref{fig:smallAllE2C1}
but with $G(\theta) = 1$ is shown in Fig.~\ref{fig:smallAllE1C1}-A.
Notice the discontinuity-induced bifurcation curve is unchanged.
This is because the difference in the two functions of $G(\theta)$ is $O(\theta^2)$
and the discontinuity-induced bifurcation is local to the origin.
The quadratic terms affect the criticality of the bifurcation;
indeed changing the function $G$ has reversed the criticality
on the bifurcation curve, Fig.~\ref{fig:smallAllE1C1}.
As above, here a curve of saddle-node bifurcations of zigzag periodic orbits
emanates from the codimension-two point at which the criticality changes.
The saddle-node bifurcation curve is shown in Fig.~\ref{fig:smallAllE1C1}-A
up until the periodic orbit at the bifurcation is no longer physically meaningful
(where it includes $\theta$-values greater than $\frac{\pi}{2}$).
For $G(\theta) = 1$ the equilibria created
at $a=1$ are stable and for this reason no
homoclinic connection forms analogous to that for $G(\theta) = \cos(\theta)$.

%---------------------------------------------------------------------
\subsection{Series expansions in $s$ and $\tau$}
\label{sub:EXPAN}

In this section we derive asymptotic expressions for bifurcations relating
to zigzag dynamics that were identified numerically in the previous section
in order to gain a greater understanding of the bifurcations.
The methodology we employ requires the period of zigzag oscillations to be small.
For this reason we assume that both $\tau$ and $s$ are small,
which is consistent with observations that zigzag dynamics occurs for small $\tau$ and $s$,
and obtain explicit expressions for
zigzag orbits as series expansions in $s$, $\tau$ and $t$.

%%%%%%%%%%%%%%%%%%%%%%%%%%%%%%%%%%%%%%%%%%%%%%%%%%%%%%%%%%%%%
\begin{figure}[b!]
\begin{center}
\setlength{\unitlength}{1cm}
\begin{picture}(8,4.3)
\put(0,0){\includegraphics[width=8cm,height=4.3cm]{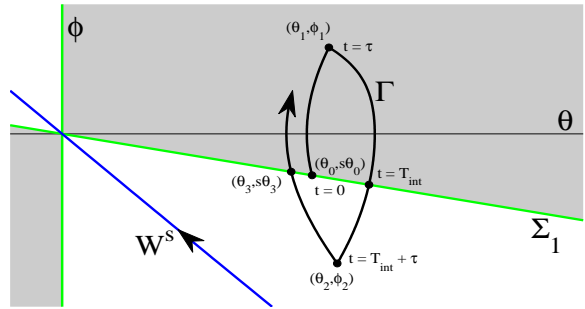}}
\end{picture}
\caption{
\label{fig:expansionSchemC1}
A sketch of the forward orbit, $\Gamma$,
of a point, $(\theta_0,s\theta_0)$, on $\Sigma_1$.
}
\end{center}
\end{figure}
%%%%%%%%%%%%%%%%%%%%%%%%%%%%%%%%%%%%%%%%%%%%%%%%%%%%%%%%%%%%%

Let $\Gamma$ be the forward orbit of a point $(\theta_0,s\theta_0)$
at $t=0$, see Fig.~\ref{fig:expansionSchemC1}.
Let $(\theta_1,\phi_1)$ denote the location of $\Gamma$ at $t = \tau$.
This is the first switching point of $\Gamma$.
Let $T_{\rm int}$ be the next time at which $\Gamma$ intersects $\Sigma_1$,
if such an intersection exists.
Then the second switching point occurs at $t = T_{\rm int} + \tau$ and
we denote this point by $(\theta_2,\phi_2)$.
For small $s$ and $\tau$ it is reasonable to assume that
the second switching point lies in the OFF region above $W^s$, so then
$\Gamma$ will exit the OFF region through $\Sigma_1$ at some point $(\theta_3,s\theta_3)$.
Naturally we are interested in the difference between $\theta_3$ and $\theta_0$
as this indicates whether $\Gamma$ is approaching or moving away from the origin.
Algebraically it is easier to instead compute the change in the
Hamiltonian, $H(\theta,\phi)$ (\ref{eq:hamiltonian}),
between the two points.
Time evolution of the OFF system does not change the Hamiltonian,
so it is equivalent to look at
\begin{equation}
\Delta H = H_2 - H_1 \;,
\label{eq:DHdef}
\end{equation}
where
\begin{equation}
H_1 = H(\theta_1,\phi_1) \;, \qquad H_2 = H(\theta_2,\phi_2) \;.
\nonumber
\end{equation}

We now compute the first few terms of $\Delta H$ as a series expansion in $s$ and $\tau$.
To derive the expansion we first directly use the
governing differential equations (\ref{eq:e})-(\ref{eq:c1})
to express $\Gamma$ as a series in $s$, $\tau$ and $t$,
then solve for $T_{\rm int}$, and finally evaluate $H_1$ and $H_2$.
Since the system under consideration includes time-delayed switching,
for our purposes this approach is preferable to
expanding $\tau$ within the differential equations
with the idea of reducing the delay differential equations
to a $\tau$-dependent ODE system.
A note with regards to notation: we use $O(|\cdot|^k)$
to denote terms of order $k$ or higher in the given variables.
We assume that $s$ and the period of oscillations are both $O(\tau)$.

The orbit $\Gamma$ is the solution to the initial value problem,
(\ref{eq:e})-(\ref{eq:c1}) with $(\theta(0),\phi(0)) = (\theta_0,s\theta_0)$.
For $t \in [0,\tau]$, $\Gamma$ is governed by the OFF system.
By substituting a series of the form
$\theta(t) = \sum_i \sum_j c_{ij}(\theta_0) s^i t^j$
into (\ref{eq:e}) with $F=0$ and solving for the coefficients we obtain
\begin{equation}
\theta(t) = \theta_0 + \theta_0 s t + \frac{1}{2} \sin(\theta_0) t^2 + O(|s,t|^4) \;,
\label{eq:expSolnOFF}
\end{equation}
which is valid for $t \in [0,\tau]$.
Substituting $t = \tau$ into (\ref{eq:expSolnOFF}) yields
\begin{eqnarray}
\theta_1 &=& \theta_0 + \theta_0 s \tau +
\frac{1}{2} \sin(\theta_0) \tau^2 + O(|s,\tau|^4) \;, \label{eq:theta1} \\
\phi_1 &=& \theta_0 s + \sin(\theta_0) \tau + O(|s,\tau|^3) \;. \label{eq:phi1}
\end{eqnarray}
Proceeding in a similar fashion we obtain an
explicit expression of $\theta(t)$ for $t > \tau$.
We provide further comments on the derivation below.
The expansions are centred about
$(\theta_1,\phi_1)$ instead of $(\theta_0,s\theta_0)$
because this provides some simplification
and then powers of $t-\tau$ naturally appear.
For small $s$, $\tau$ and $t$ and assuming $T_{\rm int} = O(\tau)$, $\Gamma$ is given by
{\small
\begin{equation}
\theta(t) = \left\{ 
\begin{array}{lc}
\theta_1 + \big( s \theta_1 + \sin(\theta_1) \tau \big) (t-\tau) +
\frac{1}{2} \sin(\theta_1) (t-\tau)^2 + O(|s,\tau,t|^4) \;, & t \in [0,\tau] \\
\theta_1 + \big( s \theta_1 + \sin(\theta_1) \tau \big) (t-\tau) +
\big( \alpha_3 + \alpha_4 s \big) (t-\tau)^2 +
\alpha_6 (t-\tau)^3 + O(|s,\tau,t|^4) \;, & t \in [\tau,2 \tau] \\
\multicolumn{2}{l}{\theta_1 + \hat{\alpha}_1 \tau^3 +
\big( s \theta_1 + \sin(\theta_1) \tau + \hat{\alpha}_2 \tau^2 \big) (t-\tau) +
\big( \alpha_3 + \alpha_4 s + \hat{\alpha}_5 \tau \big) (t-\tau)^2 +
\hat{\alpha}_6 (t-\tau)^3 + O(|s,\tau,t|^4) \;, } \\
& t \in [2 \tau,T_{\rm int} + \tau] \\
\end{array} \right.
\label{eq:expSolnAllZC1}
\end{equation}
}where
\begin{eqnarray}
\hat{\alpha}_1(\theta) &=& -\frac{1}{6} a b \theta G(\theta)^2 \;, \nonumber \\
\hat{\alpha}_2(\theta) &=& \frac{1}{2} a b \theta G(\theta)^2 \;, \nonumber \\
\alpha_3(\theta) &=& -\frac{1}{2} \Big( a \theta G(\theta) - \sin(\theta) \Big) \;, \nonumber \\
\alpha_4(\theta) &=& -\frac{1}{2} b \theta G(\theta) \;, \nonumber \\
\hat{\alpha}_5(\theta) &=& -\frac{1}{2} a b \theta G(\theta)^2 \;, \nonumber \\
\alpha_6(\theta) &=& -\frac{1}{6} b \sin(\theta) G(\theta)^2 \;, \nonumber \\
\hat{\alpha}_6(\theta) &=& 
\frac{1}{6} b G(\theta) \big( a \theta G(\theta) - \sin(\theta) \big) \;, \nonumber
\end{eqnarray}
and in (\ref{eq:expSolnAllZC1}) the $\alpha$'s are evaluated
at $\theta = \theta_1$.

The top-most expression of (\ref{eq:expSolnAllZC1}) follows from combining
(\ref{eq:expSolnOFF}) and (\ref{eq:theta1}).
The middle expression of (\ref{eq:expSolnAllZC1}) is obtained from the ON system and
substituting the top-most expression in place of $\theta(t-\tau)$
and using the time derivative of the top-most expression for $\phi(t-\tau)$.
For this reason the middle expression is valid only for $t \in [\tau,2 \tau]$.
Then using the middle expression for the time-delayed 
components, $\theta(t-\tau)$ and $\phi(t-\tau)$,
we obtain the bottom expression which is valid for $t \in [2 \tau, 3 \tau]$,
if and only if $T_{\rm int} > 2 \tau$.
Continuing in this fashion one may build up an explicit expression
for $\Gamma$ to any desired order in $s$, $\tau$ and $t$
with formulae that are valid in $[n \tau, (n+1) \tau]$, with $n \in \mathbb{Z}$,
until $t = T_{\rm int} + \tau$, beyond which the OFF system governs $\Gamma$
for some time.
The solution to the ON system, $\theta(t)$, is one degree
more differentiable than $\theta(t-\tau)$.
Consequently, $\theta(t)$ is $C^{n-1}$ at each $t=n\tau$ with
$n \in \mathbb{Z}$ and $n\tau < T_{\rm int} + \tau$.
Therefore the series solution to $\Gamma$ for $t \in [3 \tau, 4 \tau]$
differs from the bottom expression of (\ref{eq:expSolnAllZC1})
only in terms that are $O(|s,\tau,t|^4)$.
Since (\ref{eq:expSolnAllZC1}) has no explicitly stated quartic or higher order
terms the bottom expression is valid for all
$t \in [2 \tau,T_{\rm int} + \tau]$ as stated.

The intersection time, $T_{\rm int}$, is defined by
$\phi(T_{\rm int}) = s \theta(T_{\rm int})$.
From (\ref{eq:expSolnAllZC1}) we obtain
\begin{equation}
T_{\rm int} = \left\{ \begin{array}{lc}
\xi_1 \tau + \xi_2 s\tau + \xi_3 \tau^2 + O(|s,\tau|^3) \;, &
\phi(2\tau) - s \theta(2\tau) \le 0 \\
\xi_1 \tau + \xi_2 s\tau + \hat{\xi}_3 \tau^2 + O(|s,\tau|^3) \;, &
\phi(2\tau) - s \theta(2\tau) \ge 0
\end{array} \right. \;,
\label{eq:TintC1}
\end{equation}
where
\begin{eqnarray}
\xi_1 &=& \frac{a \theta_1 G(\theta_1)}{a \theta_1 G(\theta_1) - \sin(\theta_1)} \;, \nonumber \\
\xi_2 &=& \frac{-b \theta_1 \sin(\theta_1) G(\theta_1)}
{\big( a \theta_1 G(\theta_1) - \sin(\theta_1) \big)^2} \;, \nonumber \\
\xi_3 &=& \frac{-\frac{1}{2} b \sin^3(\theta_1) G(\theta_1)}
{\big( a \theta_1 G(\theta_1) - \sin(\theta_1) \big)^3} \;, \nonumber \\
\hat{\xi}_3 &=& \frac{1}{2} b G(\theta_1) \frac{
\sin^2(\theta_1) -
3 a \theta_1 \sin(\theta_1) G(\theta_1) +
a^2 \theta_1^2 G(\theta_1)^2}{\big( a \theta_1 G(\theta_1) - \sin(\theta_1) \big)^2} \;. \nonumber
\end{eqnarray}
Note that the denominators of the $\xi_i$
share the same roots as (\ref{eq:thONgrazC1}) when $s=0$.
As discussed in \S\ref{sub:NO},
a root of (\ref{eq:thONgrazC1}) corresponds to a boundary of a sliding
region of the system in the absence of delay.

From the expression for $\xi_1$ it follows that
if $a > \frac{\sin(\theta_1)}{\theta_1 G(\theta_1)}$,
i.e.~the control is sufficiently strong,
and $s$ and $\tau$ are sufficiently small,
then $\Gamma$ indeed intersects $\Sigma_1$ at a time $t = T_{\rm int} > \tau$.
If $G(\theta) = \cos(\theta)$,
this is equivalent to requiring $a > 1$ and
$\theta_0 < \theta_{\rm cos}^*$.
If $G(\theta) = 1$,
this condition is satisfied if $a > 1$ or $\theta_1^* < \theta_0 < \frac{\pi}{2}$.

$\Delta H$ is derived by evaluating $H(\theta(t),\phi(t))$
at $t = \tau$ and $t = T_{\rm int} + \tau$, and taking the difference.
The above expressions lead to
\begin{equation}
\Delta H = \left\{ \begin{array}{lc}
\zeta_1 s\tau + \zeta_2 \tau^2 + \zeta_3 s^2\tau +
\zeta_4 s \tau^2 + \zeta_5 \tau^3 + O(|s,\tau|^4) \;, &
T_{\rm int} \le 2\tau \\
\zeta_1 s\tau + \zeta_2 \tau^2 + \zeta_3 s^2\tau +
\hat{\zeta}_4 s \tau^2 + \hat{\zeta}_5 \tau^3 + O(|s,\tau|^4) \;, &
T_{\rm int} \ge 2\tau
\end{array} \right. \;,
\label{eq:DH}
\end{equation}
where
{\small
\begin{eqnarray}
\zeta_1 &=& \frac{-a^2 \theta_1^3 G(\theta_1)^2}
{a \theta_1 G(\theta_1) - \sin(\theta_1)} \;, \nonumber \\
\zeta_2 &=& -a^2 \theta_1^2 G(\theta_1)^2
\frac{\sin(\theta_1) - \frac{1}{2} a \theta_1 G(\theta_1)}
{a \theta_1 G(\theta_1) - \sin(\theta_1)} \;, \nonumber \\
\zeta_3 &=& 2 a b \theta_1^3 G(\theta_1)^2
\frac{\sin(\theta_1) - \frac{1}{2} a \theta_1 G(\theta_1)}
{\Big( a \theta_1 G(\theta_1) - \sin(\theta_1) \Big)^2} \;, \nonumber \\
\zeta_4 &=& -2 a b \theta_1^2 G(\theta_1)^2
\frac{\sin^3(\theta_1) - \frac{11}{4} a \theta_1 \sin^2(\theta_1) G(\theta_1) +
2 a^2 \theta_1^2 \sin(\theta_1) G(\theta_1)^2 -
\frac{1}{2} a^3 \theta_1^3 G(\theta_1)^3}
{\Big( a \theta_1 G(\theta_1) - \sin(\theta_1) \Big)^3} \;, \nonumber \\
\zeta_5 &=& -\frac{1}{6} a b \theta_1 \sin(\theta_1) G(\theta_1)^2
\frac{\sin^3(\theta_1) - 6 a \theta_1 \sin^2(\theta_1) G(\theta_1) +
8 a^2 \theta_1^2 \sin(\theta_1) G(\theta_1)^2 - 3 a^3 \theta_1^3 G(\theta_1)^3}
{\Big( a \theta_1 G(\theta_1) - \sin(\theta_1) \Big)^3} \;, \nonumber \\
\hat{\zeta}_4 &=& 2 a b \theta_1^2 G(\theta_1)^2
\frac{\sin^2(\theta_1) - \frac{3}{4} a \theta_1 \sin(\theta_1) G(\theta_1) +
\frac{1}{4} a^2 \theta_1^2 G(\theta_1)^2}
{\Big( a \theta_1 G(\theta_1) - \sin(\theta_1) \Big)^2} \;, \nonumber \\
\hat{\zeta}_5 &=& \frac{1}{6} a b \theta_1 G(\theta_1)^2
\frac{\sin^3(\theta_1) + 7 a \theta_1 \sin^2(\theta_1) G(\theta_1) -
9 a^2 \theta_1^2 \sin(\theta_1) G(\theta_1)^2 +
3 a^3 \theta_1^3 G(\theta_1)^3}
{\Big( a \theta_1 G(\theta_1) - \sin(\theta_1) \Big)^2} \;. \nonumber
\end{eqnarray}
}

From (\ref{eq:DH}) with $T_{\rm int} \ge 2\tau$ we may characterize
zigzag dynamics described in the previous section
($T_{\rm int} < 2\tau$ corresponds to relatively large values of $a$).
$\Gamma$ is a zigzag periodic orbit when $\Delta H = 0$.
The discontinuity-induced bifurcations of the previous section correspond to
the creation of a zigzag periodic orbit at the origin.
Therefore this bifurcation corresponds to $\Delta H = 0$
at an arbitrarily small value of $\theta_1$.
Expanding $\Delta H$ in terms of $\theta_1$
allows us to determine the nature of the discontinuity-induced bifurcations:
\begin{eqnarray}
\Delta H  &=& \bigg( -\frac{a^2}{a-1} s \tau +
\frac{a^2(a-2)}{2(a-1)} \tau^2 -
\frac{ab(a-2)}{(a-1)^2} s^2 \tau +
\frac{ab(a^2-3a+4)}{2(a-1)^2} s \tau^2 \nonumber \\
&&+~\frac{ab(3a^3-9a^2+7a+1)}{6(a-1)^2} \tau^3 + O(|s,\tau|^4) \bigg)
\theta_1^2 + O(\theta_1^4) \;.
\label{eq:DH2}
\end{eqnarray}
Note, (\ref{eq:DH2}) is
valid for both $G(\theta) = \cos(\theta)$ and $G(\theta) = 1$.
The lowest order term (i.e.~the $\theta_1^2$ term) is zero when
\begin{equation}
\tau = \frac{2}{a-2} s -
\frac{2(2+8a-15a^2+6a^3)}{3a(a-1)(a-2)^3} bs^2 + O(s^3) \;.
\label{eq:HBapprox}
\end{equation}
An omission of $O(s^3)$ terms in (\ref{eq:HBapprox}) 
gives an approximation to the occurrence of the discontinuity-induced bifurcations
and is shown in Figs.~\ref{fig:smallAllE2C1} and \ref{fig:smallAllE1C1}
(the dashed curves).
As shown in these figures the approximation agrees well with the numerical results.
The approximation was obtained from series expansions in $t$ and $\tau$
and for this reason fits the numerics less precisely
for values of $a$ near $1$, because here the period of the orbit is relatively large,
and similarly worsens with increasing $\tau$.

Since the second lowest order term in (\ref{eq:DH2}) is
of order two more than the previous term in $\theta$,
whenever this term is nonzero at a discontinuity-induced bifurcation
the bifurcating zigzag periodic orbit grows at a rate
proportional to the square-root of a non-degenerate change in 
the control parameters, $a$ and $b$.
This conclusion agrees with the bifurcation diagrams of
Figs.~\ref{fig:smallAllE2C1} and \ref{fig:smallAllE1C1}.
Furthermore the criticality of the discontinuity-induced bifurcation
is determined by the sign of the $\theta_1^4$ term.
Unlike the quadratic term, this term is dependent upon $G''(0)$
and hence differs for $G(\theta) = \cos(\theta)$ and $G(\theta) = 1$.
When $G(\theta) = \cos(\theta)$, the $\theta_1^4$ term vanishes when
{\small
\begin{equation}
\tau = \frac{3a-5}{3a^2-8a+6} s +
\frac{b \big( -243 a^6 + 1674 a^5 - 4491 a^4 + 5862 a^3 - 3611 a^2 + 642 a + 175 \big)}
{6 a (a-1) (3a^2-8a+6)^3} s^2 + O(s^3) \;.
\label{eq:critHBE2}
\end{equation}
}
Omitting $O(s^3)$ terms,
(\ref{eq:critHBE2}) intersects the discontinuity-induced bifurcation curve
(\ref{eq:HBapprox}) at the point in
Fig.~\ref{fig:smallAllE2C1}-A indicated by a circle.
Similarly when $G(\theta) = 1$, the $\theta_1^4$ term vanishes when
\begin{equation}
\tau = -\frac{2}{a} s
- \frac{2 b (a^2-2) (3a-1)}{3 a^4 (a-1)} s^2 + O(s^3) \;.
\label{eq:critHBE1}
\end{equation}

The dashed curves of Figs.~\ref{fig:smallAllE2C1} and \ref{fig:smallAllE1C1}
that approximate saddle-node bifurcations of the zigzag periodic orbits
were computed numerically using quadratic and cubic terms of (\ref{eq:DH}).

%---------------------------------------------------------------------
\subsection{Homoclinic bifurcations}
\label{sub:HOMOC}

\begin{figure}[b!]
\begin{center}
\setlength{\unitlength}{1cm}
\begin{picture}(8,4.3)
\put(0,0){\includegraphics[width=8cm,height=4.3cm]{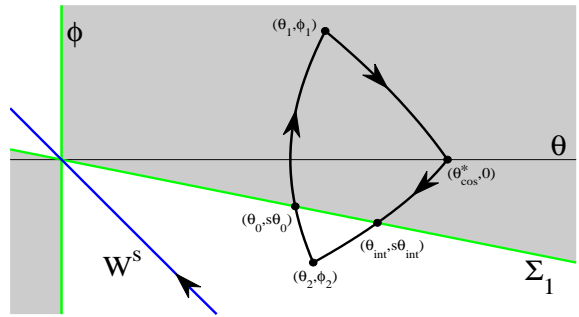}}
\end{picture}
\caption{
A sketch of a zigzag trajectory that is homoclinic to
the equilibrium, $(\theta_{\rm cos}^*,0)$.
\label{fig:homocSchem}
}
\end{center}
\end{figure}

Here we derive, to lowest order in $s$ and $\tau$, the curve of homoclinic
bifurcations shown in Fig.~\ref{fig:smallAllE2C1}.
At such a homoclinic bifurcation there exists the
homoclinic connection shown in Fig.~\ref{fig:homocSchem}.
This connection does not exist for relatively large values of $\tau$
because in that case trajectories tend to intersect $W^s$ and undergo spiral motion.

%%%%%%%%%%%%%%%%%%%%%%%%%%%%%%%%%%%%%%%%%%%%%%%%%%%%%%%%%%%%%
%%%%%%%%%%%%%%%%%%%%%%%%%%%%%%%%%%%%%%%%%%%%%%%%%%%%%%%%%%%%%

To approximate this connection we begin by deriving
approximations to the stable and unstable manifolds of $(\theta_{\rm cos}^*,0)$.
The eigenvalues and eigenvectors of the linearization
about $(\theta_{\rm cos}^*,0)$ when $\tau = 0$
are obtained by elementary calculations and give:
\begin{eqnarray}
\lambda^{\pm} &=& -\frac{b \cos(\theta_{\rm cos}^*)}{2} \pm
\sqrt{\frac{b^2 \cos^2(\theta_{\rm cos}^*)}{4} +
\frac{2\theta_{\rm cos}^* - \sin(2\theta_{\rm cos}^*)}
{2\theta_{\rm cos}^* \cos(\theta_{\rm cos}^*)}} \;,
\label{eq:sadEigVal} \\
v^{\pm} &=& \left[ \begin{array}{c} 1 \\ \lambda_1^{\pm} \end{array} \right] \;,
\label{eq:sadEigVec}
\end{eqnarray}
where $v^{\pm}$ is a vector in the $(\theta,\phi)$-plane.
The section of the orbit between the equilibrium, $(\theta_{\rm cos}^*,0)$,
and the switching point, $(\theta_2,\phi_2)$, see Fig.~\ref{fig:homocSchem},
coincides with the unstable manifold of $(\theta_{\rm cos}^*,0)$.
From (\ref{eq:sadEigVal}) and (\ref{eq:sadEigVec}),
this curve may be written as
\begin{equation}
\phi(\theta) = \lambda^+ (\theta-\theta_{\rm cos}^*) + O(|\theta-\theta_{\rm cos}^*,\tau|^2) \;.
\label{eq:sadUnstabMan}
\end{equation}

We denote the intersection of (\ref{eq:sadUnstabMan})
with $\Sigma_1$ $(\phi = s \theta)$, by $(\theta_{\rm int},s\theta_{\rm int})$
and from (\ref{eq:sadUnstabMan}) obtain
\begin{equation}
\theta_{\rm int} = \theta_{\rm cos}^* \left( 1 + \frac{s}{\lambda^+} \right) + O(s^2) \;.
\label{eq:thInt}
\end{equation}
Near $\Sigma_1$ the time-delay has more effect than near $(\theta_{\rm cos}^*,0)$.
Nevertheless, from the same series expansion methods as the previous section,
it follows that
\begin{equation}
\theta_2 = \theta_{\rm int} + O(|s,\tau|^2) \;,
\label{eq:th2HC}
\end{equation}
(i.e.~the difference between $\theta_2$ and $\theta_{\rm int}$
is negligible for this calculation).
From a series expansion in $s$ and $t$ of the solution to the OFF system
between the two switching points $(\theta_2,\phi_2)$ and $(\theta_1,\phi_1)$
we obtain, using (\ref{eq:th2HC}),
\begin{eqnarray}
\theta_1 &=& \theta_{\rm int} + O(|s,\tau|^2) \;,
\label{eq:th1HC} \\
\phi_1 &=& \theta_{\rm int} s + \sin(\theta_{\rm int}) \tau + O(|s,\tau|^2) \;.
\label{eq:phi1HC}
\end{eqnarray}
Lastly, the section of the orbit between $(\theta_1,\phi_1)$ and
$(\theta_{\rm cos}^*,0)$ is the stable manifold of the equilibrium
so given by
\begin{equation}
\phi(\theta) = \lambda_1^- (\theta-\theta_{\rm cos}^*) + O(|\theta-\theta_{\rm cos}^*|^2) \;.
\label{eq:sadStabMan}
\end{equation}
The homoclinic connection exists when the point $(\theta_1,\phi_1)$,
given by (\ref{eq:th1HC}) and (\ref{eq:phi1HC}), satisfies (\ref{eq:sadStabMan}).
Using also (\ref{eq:thInt}), we arrive at
(after some manipulation):
\begin{equation}
\tau = -\frac{1}{a} \left(
\frac{2}{\cos(\theta_{\rm cos}^*)} + \frac{b}{\lambda_1^+} \right) s + O(s^2) \;.
\label{eq:HCexistence}
\end{equation}
as a condition for the existence of the zigzag homoclinic orbit.
The approximation obtained by dropping $O(s^2)$ terms in (\ref{eq:HCexistence})
is shown in Fig.~\ref{fig:smallAllE2C1}
and matches well to the numerical results.

%=====================================================================
\section{Simple and Complex Behaviour for Large Delay}
\label{sec:LARGE}

Dynamical behaviour exhibited by the system
for small delay was described in the previous section.
However we did not provide a complete description of dynamics
in the case that $a$ is slightly larger than $1$.
In this case complex behaviour may occur
that appears to persist for larger $\tau$.
In section \S\ref{sub:BURST} we demonstrate that solutions of the system
may exhibit distinct behaviours on different time-scales
in a manner akin to bursting in neuron models.
In section \S\ref{sub:FOUR} we describe four 
bifurcations that indicate behaviour near the origin when
the delay time is large.

%---------------------------------------------------------------------
\subsection{Bursting-like dynamics}
\label{sub:BURST}

\begin{figure}[b!]
\begin{center}
\setlength{\unitlength}{1cm}
\begin{picture}(13,7.3)
\put(0,0){\includegraphics[width=13cm,height=7.3cm]{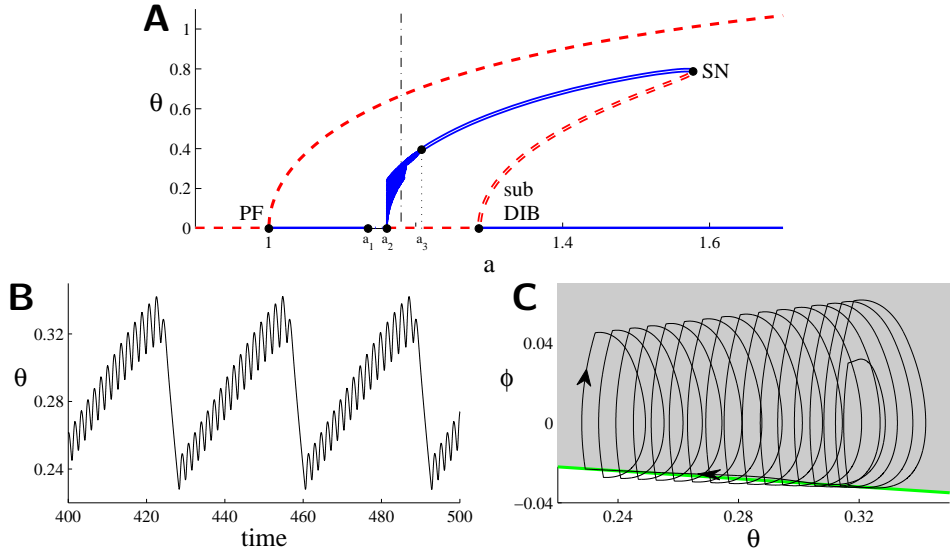}}
\put(1.9,7){\large \sf \bfseries A}
\put(.1,3.3){\large \sf \bfseries B}
\put(6.8,3.3){\large \sf \bfseries C}
\end{picture}
\caption{
Panel A is a bifurcation diagram of (\ref{eq:e})-(\ref{eq:c1})
when $\tau=0.3$, $s=-0.1$, $b=2$ and $G(\theta) = \cos(\theta)$.
The abbreviations and line styles are explained in the caption of
Fig.~\ref{fig:smallAllE2C1}-B.
Panel B shows a partial time series of an orbit with transients decayed
when $a = 1.18$ (corresponding to the dash-dot line in panel A).
Panel C shows this orbit in the $(\theta,\phi)$-plane.
For values of $a$ near $1.18$, in panel A we have indicated $\theta$-values
at which the orbit crosses $\Sigma_1$ and then immediately enters
the neighbouring ON region.
\label{fig:burst}
}
\end{center}
\end{figure}

Fig.~\ref{fig:burst}-A illustrates dynamics of (\ref{eq:e})-(\ref{eq:c1})
when $\tau = 0.3$ and $s = -0.1$ for different values of $a$.
The majority of the dynamics indicated by this plot 
match the predictions of Fig.~\ref{fig:smallAllE2C1}
for behaviour for small $\tau$.
Specifically there are stable and unstable zigzag periodic
orbits that grow in size proportional to the square root of change in $a$
and that collide and annihilate in a saddle-node bifurcation.
However for $a_2 < a < a_3$
(where $a_2 \approx 1.161$ and $a_3 \approx 1.208$, see Fig.~\ref{fig:burst}-A),
numerical simulations reveal a complicated attracting set.
This set is shown in panels B and C of Fig.~\ref{fig:burst} for $a = 1.18$.

%%%%%%%%%%%%%%%%%%%%%%%%%%%%%%%%%%%%%%%%%%%%%%%%%%%%%%%%%%%%%
%%%%%%%%%%%%%%%%%%%%%%%%%%%%%%%%%%%%%%%%%%%%%%%%%%%%%%%%%%%%%

To describe this set and its associated dynamics consider
the forward orbit of a point on $\Sigma_1$, near the origin.
For the parameter values of Fig.~\ref{fig:burst}-B, the forward orbit
zigzags away from the origin, for some time.
Over the course of each zigzag,
the orbit spends a time greater than $\tau$ within the OFF region
that decreases for each successive zigzag.
The combination of slowly outward-moving motion (with $\dot{\theta} = O(\tau)$),
and rapid zigzag motion (of frequency $O(\frac{1}{\tau})$) continues until
$\theta \approx 0.32$ at which the orbit
enters and exits the OFF region in a time less than $\tau$.
Let $t_{\rm short} < \tau$ denote this time and
let $(\theta_{\rm short},s\theta_{\rm short})$ denote the point on $\Sigma_1$
at which the orbit exits the OFF region.
Unlike what has been typical throughout this paper,
beyond $(\theta_{\rm short},s\theta_{\rm short})$
the orbit is governed by the ON system
for the next $\tau - t_{\rm short}$ time after which
it is governed by the OFF system for the next $t_{\rm short}$ time.
Numerically we observe that because this part of the orbit follows the
OFF system for less time than it would on a typical zigzag oscillation,
it does not attain such a large $\phi$-value before the control is reapplied.
Consequently the orbit may to become trapped in the ON region,
albeit only slightly above $\Sigma_1$, for some time.
This is seen in Fig.~\ref{fig:burst}-C.
During this time the orbit rapidly approaches the unstable manifold
of the saddle equilibrium, $(\theta_{\rm cos}^*,0)$.
The unstable manifold tends to the origin
but for the parameter values of Fig.~\ref{fig:burst}-B,
the manifold intersects $\Sigma_1$ at $\theta \approx 0.236$.
This $\theta$-value is sensitive to the choice of the value of $a$,
as visible in Fig.~\ref{fig:burst}-A.

Therefore, after a zigzag oscillation involving a time in the OFF region that is less
than the delay time of the system, $\tau$,
the forward orbit endures an excursion back to $\Sigma_1$
(at $\theta \approx 0.236$ in Fig.~\ref{fig:burst}-B)
followed then by outward-moving zigzag motion
and continued repetition of this procedure.
For $a = 1.18$ the attracting set is periodic
but for different values of $a$ with $a_2 < a < a_3$,
the attracting set may be chaotic.

As the value of $a$ is decreased from $a = 1.18$,
the range of $\theta$-values over which the attracting set exists, increases,
until at $a = a_2$ the unstable manifold of $(\theta_{\rm cos}^*,0)$
no longer intersects $\Sigma_1$.
For $a < a_2$ trajectories following the manifold limit directly to the origin.
It is interesting to note that for $a_1 < a < a_2$ ($a_1 \approx 1.135$),
the forward orbit of a point on $\Sigma_1$ arbitrarily close to the origin
zigzags slowly away from the origin until at
a $\theta$-value of order $1$ the orbit becomes trapped in the ON region
and approaches the origin.
Consequently over this range of $a$-values the origin
is not Lyapunov stable, and hence not asymptotically stable \cite{Gl99},
but appears to be quasi-asymptotically stable in that all points
in a neighbourhood of the origin tend to origin, eventually.

The complicated attracting sets are born out of the stable zigzag
periodic orbit in a bifurcation at $a = a_3$, Fig.~\ref{fig:burst}-A.
At the bifurcation the amount of time spent by the periodic orbit in the OFF region
is exactly the delay time, $\tau$.
This type of discontinuity-induced bifurcation was analyzed in \cite{Si06}
(see also \cite{CoDi07,SiKo10b})
for a general time-delayed, piecewise-smooth system
comprised of ordinary differential equations
on each side of the switching manifold.
In that paper it was proved that in a neighbourhood of the bifurcation
a Poincar\'{e} map is generically piecewise-smooth continuous,
and to lowest order piecewise-linear.
However, for the system studied in this paper
phase space is infinite-dimensional
so it is not clear how to define a Poincar\'{e}
section that captures all oscillatory motions local to the bifurcation.
We leave for future work a thorough investigation of
the discontinuity-induced bifurcation
characterized by a time spent in the OFF region exactly equal to
the delay time of the system.

%---------------------------------------------------------------------
\subsection{Four fundamental bifurcations for dynamics near the origin}
\label{sub:FOUR}

It is more difficult to classify dynamics of the system
when the delay time, $\tau$, is large.
For large $\tau$, spiral dynamics
(described in \S\ref{sub:BASIC}) may dominate.
Spiral dynamics cannot be analyzed by the asymptotic methods of
\S\ref{sub:EXPAN} because the time taken for a single spiral is order $1$.
Furthermore, spiral periodic orbits may undergo symmetry breaking bifurcations
followed by period-doubling cascades to complex attractors.
We leave a more complete analysis of these bifurcations for future work.
In this section we analyze the stability of the origin.

%%%%%%%%%%%%%%%%%%%%%%%%%%%%%%%%%%%%%%%%%%%%%%%%%%%%%%%%%%%%%
\begin{figure}[b!]
\begin{center}
\includegraphics[width=11cm]{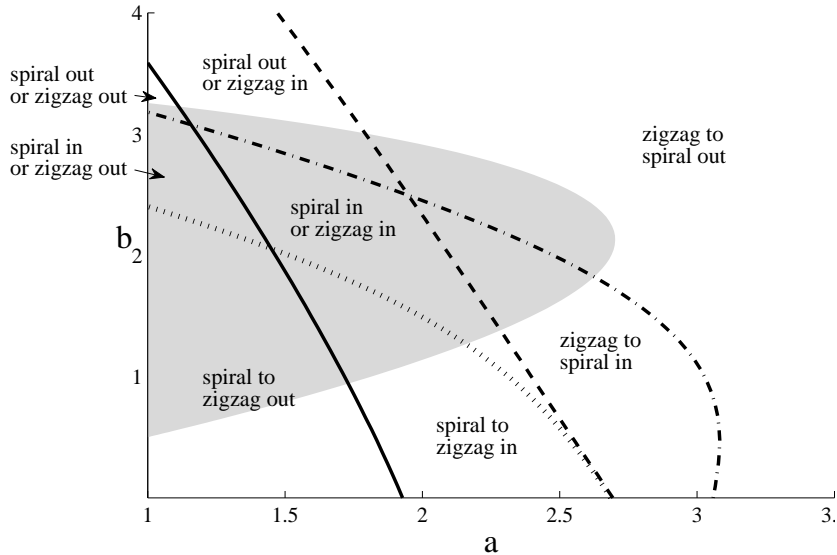}
\caption{
A bifurcation set of the linearization, (\ref{eq:linearized}),
with (\ref{eq:PD})-(\ref{eq:c1}), $\tau = 0.5$ and $s=0$.
The $(a,b)$-plane has been partitioned according to the fate
of forward evolution of the points $(1,s)$ and $(0,1)$.
By {\em zigzag in} [{\em zigzag out}] we mean that
the forward orbit of $(1,s)$ zigzags into [away from] the origin.
Similarly by {\em spiral in} [{\em spiral out}] we mean that
the forward orbit of $(0,1)$ spirals into [away from] the origin.
By {\em zigzag to spiral} we mean that the forward orbit
of $(1,s)$ undergoes spiral motion and behaves like the
forward orbit of $(0,1)$, and vice-versa for {\em spiral to zigzag}.
For instance when $(a,b) = (3,3)$ both forward orbits spiral away from the origin;
when $(a,b) = (1.5,3.5)$, the forward orbit of $(1,s)$ zigzags in to the origin
and the forward orbit of $(0,1)$ spirals away from the origin,
as in Fig.~\ref{fig:onOffc1}.
For comparison we have also indicated, by the gray D-shaped region,
the region where the origin is a stable equilibrium of the ON system.
Within this region the countable set of eigenvalues associated with the equilibrium
all have negative real part \cite{StKo00,SiKr04b}.
\label{fig:largePWLC1}
}
\end{center}
\end{figure}
%%%%%%%%%%%%%%%%%%%%%%%%%%%%%%%%%%%%%%%%%%%%%%%%%%%%%%%%%%%%%

The origin is a non-differentiable point of (\ref{eq:e})-(\ref{eq:c1})
and so does not have associated eigenvalues that determine stability.
The stability of such a point in a piecewise-smooth ODE system
is understood in two dimensions but is yet to be completely
solved in higher dimensions \cite{LiAn09,IwHa06,CaFr06,CaDe06}.
The presence of time-delay only adds complexity;
for this reason we rely on numerical simulations.
Since the sole interest here is on dynamics local to the origin,
it is sufficient to analyze the linearization of (\ref{eq:e}):
\begin{equation}
\begin{split}
\dot{\theta} &= \phi \;, \\
\dot{\phi} &= \theta + F \;,
\end{split}
\label{eq:linearized}
\end{equation}
which is valid for both $G(\theta) = \cos(\theta)$ and $G(\theta) = 1$.
Due to the scale invariance of the spatial coordinates
it suffices to look at the forward orbits
of just one point on $\Sigma_1$, say $(1,s)$,
and one point on $\Sigma_2$, say $(0,1)$.

With the goal of determining the dynamics for all combinations of the control parameters,
we perform numerical integration to study these two forward orbits.
For each, we identify, and numerically continue, two key bifurcations.
First, the forward orbit of $(1,s)$ may return to this point upon one zigzag.
In the $(a,b)$-plane this occurs along the solid curve in Fig.~\ref{fig:largePWLC1},
for which $\tau = 0.5$ and $s = 0$.
We have found that different values of $\tau$ and $s$
produce qualitatively similar pictures.
Locally, for values of $a$ and $b$ to the left of this curve,
the forward orbit of $(1,s)$ zigzags away from the origin;
to the right of the curve the orbit zigzags into the origin.
For the nonlinear system, (\ref{eq:e}),
this is the discontinuity-induced bifurcation
identified in \S\ref{sub:BIFSET} at which
two symmetric zigzag periodic orbits are created.
Second, the forward orbit of $(1,s)$ may become coincident to $W^s$
(the stable manifold of the origin for the OFF system, Fig.~\ref{fig:onOffc1}).
This occurs along the dashed curve in Fig.~\ref{fig:largePWLC1}.
For values of $a$ and $b$ to the left of this curve the forward orbit zigzags,
to the right of the curve it spirals.

For the initial point $(0,1)$
there are two bifurcations analogous to those just discussed.
Along the dash-dot curve in Fig.~\ref{fig:largePWLC1}
the forward orbit of $(1,s)$ returns to $(1,s)$ after one spiral,
along the dotted curve the orbit falls onto $W^s$
and limits upon the origin without again crossing
either of the switching manifolds.
The dash-dot curve is a bifurcation of (\ref{eq:e})
at which a symmetric spiral periodic orbit is born
in a manner akin to a Hopf bifurcation.

Note that one may choose the control parameters such that zigzag orbits approach 
the origin and spiral orbits head away from the origin
(as in Fig.~\ref{fig:onOffc1}) and vice-versa.
We have not been able to identify
an intersection between the dashed and dotted curves of Fig.~\ref{fig:largePWLC1}
for any $\tau$ and $s$, nor have been able to show
that such an intersection cannot occur.
Such an intersection could permit
for the existence of an orbit that repeatedly switches
between zigzag and spiral motion.

%=====================================================================
\section{Dynamics with the Switching Condition (\ref{eq:c2})}
\label{sec:OO}	% Origin Omitted

\begin{figure}[b!]
\begin{center}
\includegraphics[width=7.2cm,height=6cm]{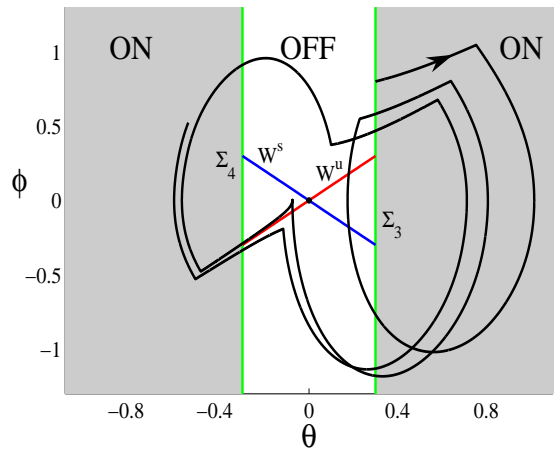}
\caption{
The $(\theta,\phi)$-plane for (\ref{eq:e})-(\ref{eq:PD}) with (\ref{eq:c2})
and $\sigma = 0.3$, $\tau = 0.5$,
$(a,b) = (2.5,4)$ and $G(\theta) = \cos(\theta)$.
A trajectory approaches an asymmetric periodic orbit.
$W^s$ and $W^u$ are the stable and unstable manifolds
of the origin for the OFF system, respectively.
\label{fig:onOffc2}
}
\end{center}
\end{figure}

In this section we analyze the system (\ref{eq:e})-(\ref{eq:PD})
with the alternative switching rule (\ref{eq:c2}), Fig.~\ref{fig:onOffc2}.
With this rule the control is removed
when the time-delayed position is near vertical
and is motivated from observations of human balancing
tasks and postural sway \cite{MiOh09}.
Switching control off near the origin may lessen ``over-control''
but eliminates the possibility of a stable vertical position.
Previous investigations have used equations of motion that are
linear in $\theta$ \cite{MiOh09,KoSt00,KoGl10}.
% [MiOh09] 1d model, brief discussion
% [KoSt00] linear eqns
% [KoGl10] considers same system but linearized in $\theta$,
% 	also no time-delay in switching rule
For typical practical applications this is justified
because the relevant range of $\theta$ values is sufficiently small.
We have chosen not to linearize in $\theta$
since our asymptotic methods do not rely on it and
so that we may study the influence of additional equilibria.

%%%%%%%%%%%%%%%%%%%%%%%%%%%%%%%%%%%%%%%%%%%%%%%%%%%%%%%%%%%%%%%%%%%%%%
%%%%%%%%%%%%%%%%%%%%%%%%%%%%%%%%%%%%%%%%%%%%%%%%%%%%%%%%%%%%%%%%%%%%%%

The system (\ref{eq:e})-(\ref{eq:PD}) with (\ref{eq:c2}) is infinite-dimensional,
so, as above, we look only at the forward orbits of points on switching manifolds.
We believe that an analysis of the fate of these orbits
provides a good understanding of the important properties of the system
for a wide range of parameter values.
Orbits that cross the OFF region in a time less than $\tau$
occur readily if $\tau$ is large relative to the width of the OFF region.
In this case typical stable dynamics are oscillations that
involve a relatively large range of $\theta$-values
which spend only a small fraction of time in the OFF region and are therefore
often similar to dynamics of simply the ON system.

%%%%%%%%%%%%%%%%%%%%%%%%%%%%%%%%%%%%%%%%%%%%%%%%%%%%%%%%%%%%%%%%%%%%%%
\begin{figure}[b!]
\begin{center}
\setlength{\unitlength}{1cm}
\begin{picture}(10,14.1)
\put(0,7.3){\includegraphics[width=11cm,height=6.8cm]{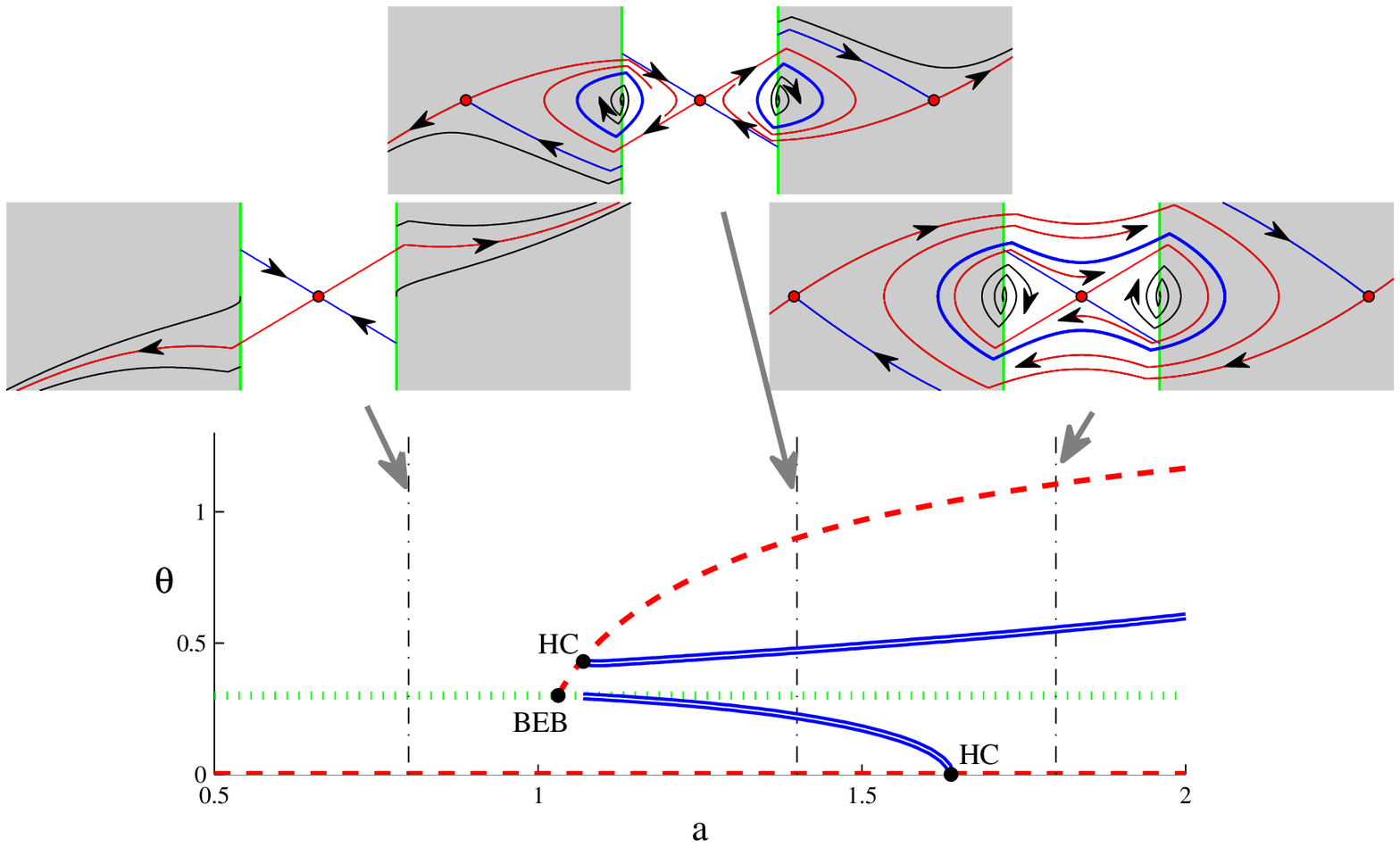}}
\put(0,0){\includegraphics[width=11cm,height=6.8cm]{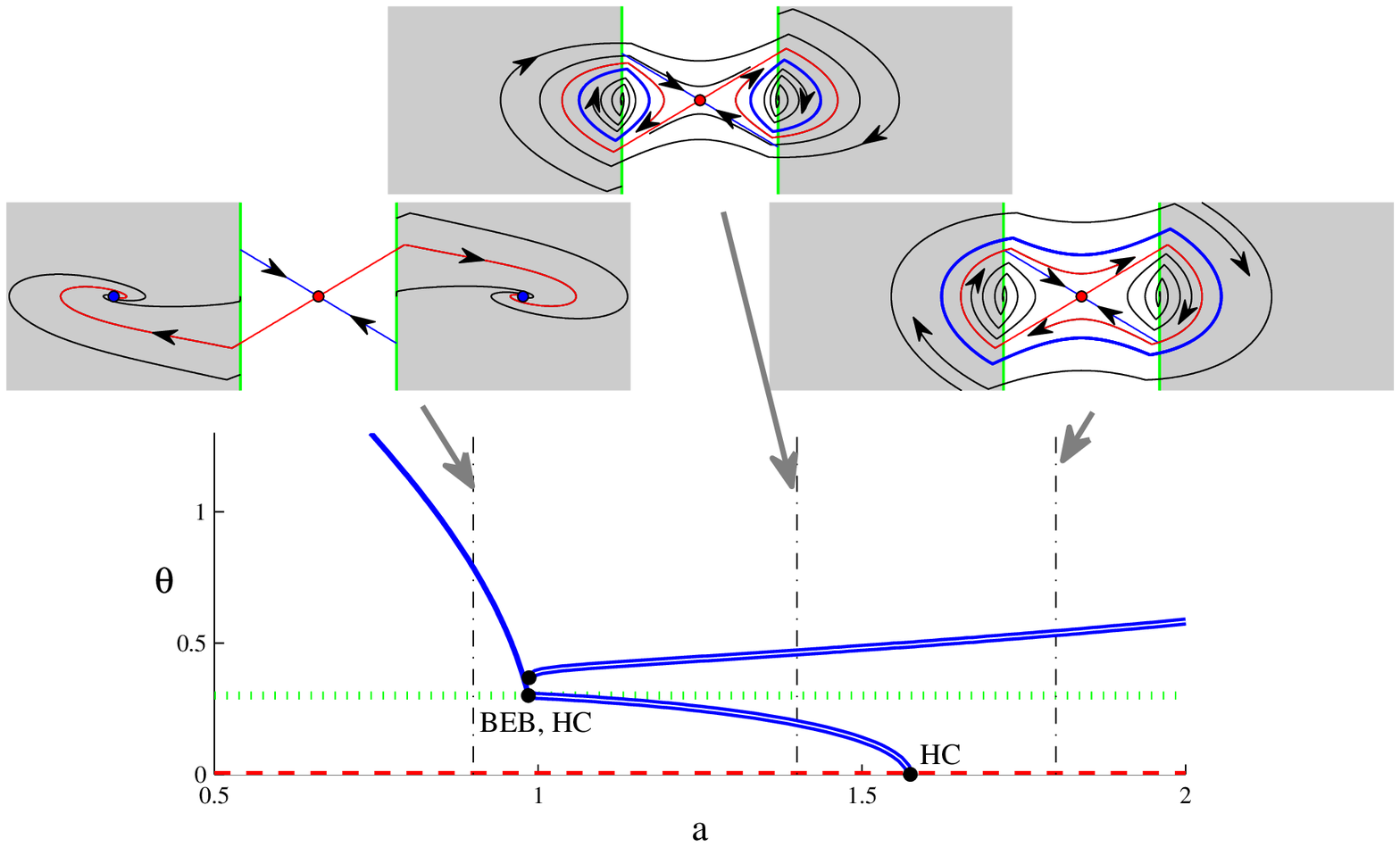}}
\put(.7,13.7){\large \sf \bfseries A}
\put(.7,6.4){\large \sf \bfseries B}
\end{picture}
\caption{
Bifurcation diagrams and representative sketches of dynamics in the
$(\theta,\phi)$-plane of the system (\ref{eq:e})-(\ref{eq:PD})
with (\ref{eq:c2}) for $\tau = 0.1$, $\sigma = 0.3$, $b = 0.5$ and
$G(\theta) = \cos(\theta)$ in panel A and $G(\theta) = 1$ in panel B.
HC - homoclinic bifurcation;
BEB - boundary equilibrium bifurcation.
Due to the symmetry of the system, dynamics for $\theta < 0$ are identical
to those for $\theta > 0$ and for this reason are not shown.
Between the homoclinic bifurcations the 
double curves indicate the maximum and minimum $\theta$-values
of a stable periodic orbit.
To the right of the right-most homoclinic bifurcation
there exists one symmetric periodic orbit;
its minimum value is not visible in the bifurcation diagrams.
Solid [dashed] curves correspond to stable [unstable] equilibria.
The dotted lines represent the switching manifold, $\theta = \sigma$.
\label{fig:smallAllC2}
}
\end{center}
\end{figure}
%%%%%%%%%%%%%%%%%%%%%%%%%%%%%%%%%%%%%%%%%%%%%%%%%%%%%%%%%%%%%%%%%%%%%%

For ease of explanation we discuss dynamics only for $\theta > 0$;
by symmetry identical dynamics occur for $\theta < 0$.
On the switching manifold $\Sigma_3$ (\ref{eq:swManC2}), $\dot{\theta} = \phi$,
thus trajectories that cross $\Sigma_3$ at, say, $(\sigma,\phi_0)$,
next enter the neighbouring ON region if $\phi_0 > 0$ and
next enter the OFF region if $\phi_0 < 0$.
For this reason we study the forward orbits of points
$(\sigma,\phi_0)$ with $\phi_0 > 0$
in view of our earlier discussion regarding initial conditions, \S\ref{sub:SWITCH}.

Fig.~\ref{fig:smallAllC2} shows bifurcation diagrams
when $\tau = 0.1$, $\sigma = 0.3$ and $b = 0.5$.
Note that the value of $\sigma$ used in this illustration is significantly
larger than values suitable for traditional human balancing tasks
so may be more meaningful in regards to mechanical applications.
However, the purpose of Fig.~\ref{fig:smallAllC2} is to highlight the
basic bifurcation structure of the model.
We have identified qualitatively similar dynamics over a wide range of parameter values,
in particular with smaller values of $\sigma$ and $\tau$.
The most parameter sensitive component of the bifurcation structure is
the value of $a$ of the right-most homoclinic bifurcation which
decreases with an increase in $\tau$.

The equilibrium of the ON system,
$(\theta_{\rm cos}^*,0)$ or $(\theta_1^*,0)$, is admissible
if and only if it lies to right of $\Sigma_3$.
Upon variation of $a$, the equilibrium collides with the switching manifold at
the point $(\sigma,0)$.
(This is referred to as a {\em boundary equilibrium bifurcation} \cite{DiBu08}.)
For values of $a$ between the two homoclinic bifurcations, Fig.~\ref{fig:smallAllC2},
there exists a stable periodic orbit
encircling $(\sigma,0)$ (and a symmetric orbit encircling $(-\sigma,0)$).
As the value of $a$ is decreased the periodic orbit is destroyed in a
homoclinic bifurcation with $(\theta_{\rm cos}^*,0)$ or $(\theta_1^*,0)$
(depending on the function $G(\theta)$).
In the case $G(\theta) = 1$ this bifurcation is coincident with
the boundary equilibrium bifurcation.
As the value of $a$ is increased the two symmetric periodic orbits
connect at the origin in a homoclinic bifurcation beyond which
there exists one stable symmetric periodic orbit
encircling the origin.
With a further increase in $a$,
or an increase in $\tau$,
numerically we have observed that this periodic orbit undergoes symmetry breaking and
period doubling to an aperiodic attractor in a fashion similar
to the system with continuous control \cite{SiKr04}.

For the remainder of this section we perform an asymptotic analysis,
along the same lines as in \S\ref{sub:EXPAN} for (\ref{eq:c1}),
to analytically derive the small amplitude periodic orbit encircling $(\sigma,0)$
for small $\tau$, and determine the rate at which the periodic orbit
grows in size with $\tau$.

%%%%%%%%%%%%%%%%%%%%%%%%%%%%%%%%%%%%%%%%%%%%%%%%%%%%%%%%%%%%%
\begin{figure}[b!]
\begin{center}
\includegraphics[width=8cm]{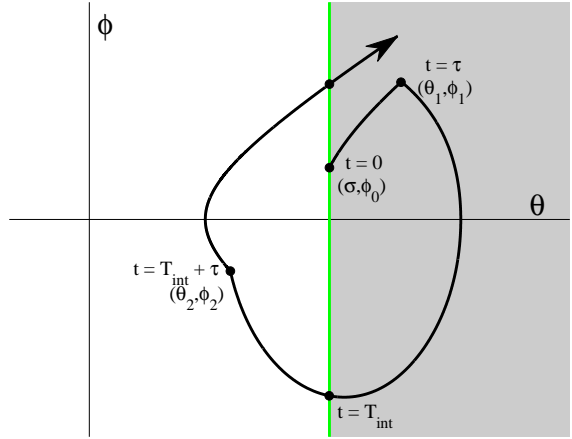}
\caption{
A sketch of the forward orbit of a point on $\Sigma_3$
for the system (\ref{eq:e})-(\ref{eq:PD}) with the switching rule (\ref{eq:c2}).
\label{fig:expansionSchemC2}
}
\end{center}
\end{figure}
%%%%%%%%%%%%%%%%%%%%%%%%%%%%%%%%%%%%%%%%%%%%%%%%%%%%%%%%%%%%%

For small $\phi_0 > 0$, let $\Gamma$ be the forward orbit of $(\sigma,\phi_0)$ at $t=0$, 
for (\ref{eq:e})-(\ref{eq:PD}) with (\ref{eq:c2})
and, recalling the discussion in \S\ref{sub:SWITCH},
assume that $\Gamma$ lies in the OFF region for all $t \in [-\tau,0]$.
Switching occurs within the ON region at $t = \tau$ and, assuming $\Gamma$
re-enters the OFF region at some later time $T_{\rm int}$, a second switching occurs at
$t = T_{\rm int} + \tau$ before $\Gamma$ exits the OFF region.
We let $(\theta_1,\phi_1)$ and $(\theta_2,\phi_2)$ denote
the respective switching points, Fig.~\ref{fig:expansionSchemC2}.
Our goal is to calculate the change in the Hamiltonian (\ref{eq:hamiltonian})
between the two switching points in order to identify periodic orbits.
Let
\begin{equation}
H_1 = H(\theta_1,\phi_1) \;, \qquad H_2 = H(\theta_2,\phi_2) \;.
\nonumber
\end{equation}
$\Gamma$ is periodic if $\Delta H \equiv H_2 - H_1 = 0$.
Since $H_1 = H(\sigma,\phi_0)$, we have
\begin{equation}
H_1 = \frac{1}{2} \phi_0^2 + \cos(\sigma) \;.
\label{eq:H1C2}
\end{equation}

As in \S\ref{sub:EXPAN},
we obtain a useful description of $\Gamma$ by substituting a series
representation of $\theta(t)$ expanded
in $\phi_0$ and $t$ into the equations of motion
and solving for the unknown coefficients:
{\small
\begin{equation}
\theta(t) = \left\{ \begin{array}{lc}
\sigma + \phi_0 t + \frac{1}{2} \sin(\sigma) t^2 + O(|\phi_0,t|^4) \;, &
t \in [0,\tau] \\
( \sigma + \phi_0 \tau + \frac{1}{2} \sin(\sigma) \tau^2 ) +
( \phi_0 + \sin(\sigma) \tau ) (t-\tau) \\
\quad+~( \alpha_3 + \frac{\alpha_4}{\sigma} \phi_0 ) (t-\tau)^2 +
\frac{\alpha_6}{G(\sigma)} (t-\tau)^3 + O(|\phi_0,\tau,t|^4) \;, &
t \in [\tau,2\tau] \\
( \sigma + \phi_0 \tau + \frac{1}{2} \sin(\sigma) \tau^2 + \hat{\alpha}_1 \tau^3 ) +
( \phi_0 + \sin(\sigma) \tau + \hat{\alpha}_2 \tau^2 ) (t-\tau) \\
\quad+~( \alpha_3 + \frac{\alpha_4}{\sigma} \phi_0 + \hat{\alpha}_5 \tau ) (t-\tau)^2 +
\hat{\alpha}_6 (t-\tau)^3 + O(|\phi_0,\tau,t|^4) \;, &
t \in [2\tau,T_{\rm int} + \tau]
\end{array} \right. \;,
\label{eq:expSolnAllZC2}
\end{equation}
}where the $\alpha$'s, listed in \S\ref{sub:EXPAN}, are evaluated at $\theta = \sigma$.
To determine $T_{\rm int}$ from $\theta(T_{\rm int}) = \sigma$ it is necessary to
consider $\phi_0 = O(\tau^{\frac{1}{2}})$ and write
\begin{equation}
T_{\rm int} = \chi_1 \phi_0 + \chi_2 \phi_0^2 + \chi_3 \tau
+ O(|\phi_0,\tau^{\frac{1}{2}}|^3) \;.
\label{eq:TintC2}
\end{equation}
Analogous to \S\ref{sub:EXPAN}, the unknown coefficients
are calculated by substituting (\ref{eq:TintC2}) into (\ref{eq:expSolnAllZC2}):
\begin{eqnarray}
\chi_1 &=& \frac{2}{a \sigma G(\sigma) - \sin(\sigma)} \;, \nonumber \\
\chi_2 &=& \frac{-\frac{2}{3} b G(\sigma)}{(a \sigma G(\sigma) - \sin(\sigma))^2} \;, \nonumber \\ 
\chi_3 &=& \frac{2 a \sigma G(\sigma)}{a \sigma G(\sigma) - \sin(\sigma)} \;, \nonumber
\end{eqnarray}
assuming $T_{\rm int} > 2 \tau$.
Notice we must have $\sin(\sigma) - a \sigma G(\sigma) < 0$
because we require $T_{\rm int} > 0$.
An evaluation of (\ref{eq:hamiltonian}) at $t = T_{\rm int} + \tau$ using
(\ref{eq:H1C2}), (\ref{eq:expSolnAllZC2}) and (\ref{eq:TintC2}) produces
\begin{equation}
\Delta H = 2 a \sigma G(\sigma) \phi_0 \tau -
\frac{\frac{2}{3} b G(\sigma)}{a \sigma G(\sigma) - \sin(\sigma)} \phi_0^3 +
O(|\phi_0,\tau^{\frac{1}{2}}|^4) \;.
\end{equation}
Consequently $\Delta H = 0$ when
\begin{equation}
\tau = \frac{b}{3 a \sigma ( a \sigma G(\sigma) - \sin(\sigma) )} \phi_0^2 + O(\phi_0^3) \;,
\label{eq:DHis0}
\end{equation}
which is consistent with our assumption that $\phi_0 = O(\tau^{\frac{1}{2}})$.

In summary, when $\Delta H = 0$, $\Gamma$ is
a small amplitude periodic orbit encircling $(\sigma,0)$,
see for instance Fig.~\ref{fig:smallAllC2}.
As the value of $\tau$ is increased from zero,
we deduce from (\ref{eq:DHis0})
this periodic orbit grows out of the point $(\sigma,0)$
with an amplitude asymptotically proportional to $\tau^{\frac{1}{2}}$.
Furthermore, at this periodic orbit
$\frac{\partial \Delta H}{\partial \phi_0} =
\frac{-\frac{4}{3} b G(\sigma)}{a \sigma G(\sigma) - \sin(\sigma)} \phi_0^2 + O(\phi_0^3)$,
which is negative-valued because $a \sigma G(\sigma) - \sin(\sigma) > 0$.
Consequently the periodic orbit is stable
matching the numerical results of Fig.~\ref{fig:smallAllC2}.

%=====================================================================
\section{Discussion}
\label{sec:CONC}

In this paper we have identified bifurcations and dynamics
in a prototypical balancing model describing planar motion of an
inverted pendulum with control that is qualitatively affected by the
combination of time-delay, discontinuity in the control, and nonlinearity.
Time-delay is fundamental to a variety of balancing problems.
It represents the time-lag between when variables are measured
and corrective forces are applied, which in human balancing tasks
typically represents neural transmission time.
Switching in the method of control has been proposed to reflect
observations of intermittent muscle movements,
to procure simplicity in mechanical systems,
or to provide a stabilizing mechanism
particularly when the time-delay is long.
Finally, terms in the equations of motion that are nonlinear in the
angle of displacement from vertical, $\theta$,
are important when the value of $\theta$ is not restricted to small values.

Previous work uses mathematical methods to analyze systems
that lack any one of these three features.
The bifurcation theory and methods of piecewise-smooth systems \cite{DiBu08}
apply to systems without time-delay.
Centre manifold reductions may be applied to models
that lack a switching condition and are smooth \cite{Ca09,SiKr04b}.
Systems that are spatially scale-invariant have been considered
in the context of balancing \cite{AsTa09}
and in general \cite{LiAn09,IwHa06,CaFr06,CaDe06};
numerical simulations are often essential in this case.
Even though (\ref{eq:e})-(\ref{eq:PD}) with either (\ref{eq:c1}) or (\ref{eq:c2})
exhibits all three of the above features
we have been able to obtain some formal results.
One simplifying aspect is that the system does not switch between two DDEs,
as in for instance \cite{Si06,SiKo10b},
but rather switches between a DDE (the ON system) and an ODE (the OFF system).
As a result, whenever an orbit spends a continuous length of time equal to or greater
than the delay time, $\tau$, governed purely by the OFF system,
its future evolution becomes independent of its location at any earlier time.
Consequently initial conditions of the system
may be thought of as points in the $(\theta,\phi)$-plane.
More specifically, since the dynamics of the OFF system is lucid,
for initial conditions we use points on the boundaries of the ON/OFF regions
at which the vector field of the OFF system points
into the neighbouring ON region, \S\ref{sub:SWITCH}.

The presence of time-delay in the switching rules
induces what we have referred to as zigzag motion.
This motion is characterized by a rapid on/off switching of the control
and corresponds to a restriction of the pendulum to one side of the vertical position.
Since zigzag motion occurs on an $O(\tau)$ time-scale,
it succumbs to the asymptotic approach,
based on piecewise Taylor series in $t$ and $\tau$,
for example (\ref{eq:expSolnAllZC1}).
We performed the asymptotic analysis
for the particular state-dependent switching rule, (\ref{eq:c1}),
where we expanded also in the switching parameter $s$.
When $-s > \tau$, and $\tau$ is not too large,
zigzag motion of the pendulum approaches the vertical position
on a relatively long time-scale, \S\ref{sub:ZIGZAG}.
This manner of stabilizing the vertical position is not possible
without a switching rule like (\ref{eq:c1}).

Nonlinearity in $\theta$ in (\ref{eq:e})-(\ref{eq:PD}) with (\ref{eq:c1}) permits
non-equilibrium, asymptotically stable invariant sets.
Using the series expansions mentioned above, we have been able to identify
periodic orbits of period $O(\tau)$,
and derive equations in terms of the system parameters that
correspond to bifurcations of these periodic orbits.
For instance zigzag periodic orbits bifurcate from the vertical
position in symmetric discontinuity-induced bifurcations.
We have analyzed the stability of these periodic orbits
both numerically, \S\ref{sub:BIFSET},
and through asymptotic expansions, \S\ref{sub:EXPAN}.
A homoclinic bifurcation was identified for small $\tau$ and $s$
in a similar fashion, \S\ref{sub:HOMOC}.
We also described a complicated bursting-like attractor in \S\ref{sub:BURST}.

For relatively large values of $\tau$ the model predicts
the pendulum to typically prefer oscillations about the vertical position.
We have referred to such motion as spiral motion
due to the nature of corresponding trajectories in the $(\theta,\phi)$-plane.
In \S\ref{sub:FOUR} we have investigated  the spiral motion numerically in
the context of the stability of the vertical position.
Since the spiral behaviour operates on long time-scales it cannot be
analyzed by the asymptotic approach described above.
In particular we found that the vertical position may be semi-stable 
in that for a fixed choice of control parameters,
spiral motion may approach the vertical position
whereas zigzag motion heads away from this position, or vice-versa.
It is interesting to note that dynamics local to the vertical position
are explained by a global analysis of the piecewise-linear system, (\ref{eq:linearized}).
This circumstance of global dynamics governing local behaviour is a common
occurrence in piecewise-smooth systems, see for instance \cite{FrPo98}.

In \S\ref{sec:OO} we studied the model with the
switching rule (\ref{eq:c2}) that turns off the control
when the controller interprets the magnitude of $\theta$ to be
less than a threshold value, $\sigma$.
In this setup the vertical position is always unstable.
Using the same asymptotic methods
we showed that as the value of $\tau$ is increased from zero
a stable periodic orbit emanates from $(\theta,\phi) = (\sigma,0)$
with an amplitude asymptotically proportional to $\tau^{\frac{1}{2}}$.
This periodic orbit also corresponds to small periodic fluctuations
of the pendulum on one side of the vertical position
due to an intermittent application of the control.
By symmetry there exists an identical stable periodic orbit
on the other side of the vertical position.
As the value of $\tau$ is increased,
typically the two periodic orbits collide in a homoclinic bifurcation
with the vertical position beyond which there exists
one symmetric stable periodic orbit corresponding to oscillations about
the vertical position.
With a further increase in $\tau$
dynamics exhibited by the system are similar to
the case where the control is constantly applied.

A future project is that of
an investigation into the effect of noise in models of the type studied here.
Some steps in this direction have already been achieved \cite{AsTa09,EuMi96,MiCa08}.
Noise may result from discrepancies in measurements of controller
or from fluctuations in muscle response
and may induce a flip-flop motion between coexisting stable solutions
or possibly have a stabilizing effect \cite{MiCa08,Ca05}.

%\bibliographystyle{unsrt}
%\bibliography{../Balance,../../PhDResearch/PhDThesis}

\end{document}